\input amstex
\documentstyle{amsppt}
\magnification = 1100
\mathsurround = 1 pt
\loadbold
\define\si{\smallskip\noindent}
\define\bi{\bigskip\noindent}
\define\pr{\text{P}}

\define\e{\varepsilon}
\topmatter
\title
Tight Markov chains and random compositions
\endtitle
\rightheadtext{Tight chains}
\author
Boris Pittel
\endauthor
\affil
Ohio State University
\endaffil
\email
bgp\@math.ohio-state.edu
\endemail
\thanks
Research supported in part by the NSF Grant DMS-0805996
\endthanks
\keywords
Random compositions, Markov chains, rarity, exponentiality, extreme
values, limit theorems
\endkeywords
\subjclass
05A15, 05A17, 11P99, 60C05, 60F05, 60J05
\endsubjclass
\address
Department of Mathematics, Ohio State University, Columbus, Ohio,  43210, USA
\endaddress
\abstract
For an ergodic  Markov chain $\{X(t)\}$ on $\Bbb N$, with a stationary 
distribution $\pi$, let $T_n>0$ denote a hitting time for $[n]^c$,
and let $X_n=X(T_n)$. Around 2005 Guy Louchard popularized a conjecture that, for $n\to
\infty$, $T_n$ is almost Geometric($p$), $p=\pi([n]^c)$, $X_n$ is almost 
stationarily distributed on $[n]^c$, and that $X_n$ and $T_n$ are almost independent,
if $p(n):=\sup_ip(i,[n]^c)\to 0$ exponentially fast. For the chains with $p(n)
\to 0$ however slowly, and with $\sup_{i,j}\,\|p(i,\cdot)-p(j,\cdot)\|_{TV}<1$,
we show that Louchard's conjecture is indeed true even for the hits of an arbitrary 
$S_n\subset\Bbb N$ with $\pi(S_n)\to 0$. More precisely,
a sequence of $k$ consecutive hit locations paired with the time elapsed since
a previous hit (for the first hit, since the starting moment)  
is approximated,
within a  total variation distance of order $k\,\sup_ip(i,S_n)$, by 
a $k$-long sequence of independent copies of $(\ell_n,t_n)$, where $\ell_n=
\text{Geometric}\,(\pi(S_n))$, $t_n$ is distributed stationarily on $S_n$, and
$\ell_n$ is independent of $t_n$.
The two conditions are easily met by the Markov chains that arose in Louchard's
studies as likely sharp approximations of two random compositions of a large
integer $\nu$, a column-convex animal (cca) composition and a Carlitz (C)
composition. We show that this approximation is indeed very sharp for each of the random
compositions, read from left to right, for as long as the sum of the remaining 
parts stays above $\ln^2 \nu$.  
Combining the two approximations, 
a composition -- by its chain, and, for $S_n=[n]^c$, the sequence of hit locations 
paired each
with a time elapsed from the previous hit 
 -- by the independent copies of $(\ell_n,t_n)$, enables us to determine the
limiting distributions of $\mu=o(\ln\nu)$ and $\mu=o(\nu^{1/2})$ 
largest parts of the random cca composition and the random C-composition,
respectively. (Submitted to Annals of Probability in August,  2009.)
\endabstract
\endtopmatter  
\document
\bi
{\bf 1. Introduction.\/} Consider a Markov chain $X(t)$ on $\Bbb N$.
Given $S\subset \Bbb N$, let $T(S)$ be the hitting time, i.e.
$T(S)=\min\{t>0\,:\,X(t)\in S\,\}$. Keilson [14] proved that if
a state $i$ is positive-recurrent, and a nested sequence $S_1\supseteq S_2\supseteq
\cdots$ is such that $i\notin S_1$ and $E_i[T(S_n)]\to\infty$, then
$$
P_i\left\{\frac{T(S_n)}{E_i[T(S_n)]}\ge t\right\}\to e^{-t}, 
\quad\forall\,t\ge 0.\tag 1.1
$$
The basic idea of the proof was that the probability of hitting $S_n$
between two consecutive returns to $i$ is small, of order $1/E_i[T(S_n)]$,
and so $T(S_n)$ is roughly the sum of the geometrically distributed
number of i.i.d.~ times between those returns to $i$. 

If a chain is ergodic, with a stationary distribution $\pi$, the condition
$E_i[T(S_n)]\to\infty$ is  met if  (and only if) 
$\pi(S_n):=\sum_{i\in S_n}\pi(i)\to 0$.
Indeed, by Derman's theorem [9] (see Durrett [10], Ch. 5), the expected number
of visits to $S_n$ between two returns to $i$ is $\pi(S_n)/\pi(i)$. So
the probability of hitting $S_n$ between two returns to $i$ is $\pi(S_n)/\pi(i)$
at most, whence $E_i[T(S_n)]\ge \pi(i)/\pi(S_n)$.
\si
 
Aldous [1] estimated accuracy of the exponential approximation of
the hitting time for a finite-state ergodic Markov chain, when an
initial state is chosen at random, in accordance with the stationary
distribution $\pi$. Roughly, the discrepancy is small if the expected
hitting time far exceeds a relaxation time $\tau=\max_i\min\{t\,:\,
\|p^t(i,\cdot)-\pi(\cdot)\|_{TV}\le\rho\}$, $\rho<1/2$. $\tau$ 
``measures the time taken for the
chain to approach stationarity'' in a sense that
$\max_i\|p^t(i,\cdot)-\pi(\cdot)\|_{TV}\le (2\rho)^{\lfloor t/\tau\rfloor}$.

Precisely because these results are so strikingly general, more
subtle questions  remain open. Is there a
geometrically distributed random variable close to $T(S)$ in terms
of the total variation distance? What is, asymptotically, the joint
distribution of the hitting time $T(S)$ and the hit location
$X(T(S))$? Is there an explicit convergence rate in terms
of the total variation distance? Are $X(T(S))$ and $T(S)$ almost independent?
How does one describe asymptotic behavior of the first $k$ visits to
the rare set $S$, if $k=k(S)$ is not too large?

For an ergodic  Markov chain $\{X(t)\}$ on $\Bbb N$, with a stationary 
distribution $\pi$, let $T_n>0$ denote a hitting time for $[n]^c=\Bbb N\setminus [n]$,
and let $X_n=X(T_n)$. Around 2005 Guy Louchard [18] popularized 
the following conjecture.
If $p(n):=\sup_ip(i,[n]^c)=O(q^n)$, $q<1$, then 
$T_n$ is almost Geometric($p$), ($p=\pi([n]^c)$), $X_n$ is almost 
stationarily distributed on $[n]^c$, and $X_n$ and $T_n$ are almost independent.
The Markov chains with $p(n)=O(q^n)$ arose in the studies of
two random compositions, Louchard [19], [20] and Louchard,
Prodinger [21] as possibly sharp approximations of those random compositions.  
Louchard's thought-provoking idea was that if the conjecture and approximability 
of each random
compositions by a chain would be proved, 
potentially  one could obtain the limiting
distributions, marginal and joint, of extreme-valued parts and, possibly,  
of other related characteristics of the random compositions.

In this paper we introduce a class of Markov chains that contains the chains
from [19]-[21] for which we can give full answers to the questions 
posed above and, in particular, fully confirm Guy Louchard's conjecture. 
We also prove that the chains in [19]-[21] indeed provide a good 
approximation of the random
compositions. The two approximations made in tandem  lead to the asymptotic distributions of 
the extreme-valued parts of the compositions, together with the convergence rates. 
\si

Let us give a more specific description of our results.
 \si
 
{\bf Definition 1.1.\/} An ergodic Markov chain on $\Bbb N$, with a transition
probability matrix $P=\{p(i,k)\}_{i,k\in\Bbb N}$ and a stationary distribution
$\pi$, is called tight
if the family of row probability measures $\{p(i,\cdot)\}_{i\in \Bbb N}$ is tight,
i.e.
$$
\lim_{n\to\infty}\sup_i\sum_{k>n}p(i,k)=0.\tag 1.2
$$
For a tight $P$, we will prove that if $\emptyset\neq S_n\subset \Bbb N$ is such that $\pi(S_n)\to 0$, then
uniformly for all initial states $i$,
$$
E_i[T^k(S_n)]\sim \frac{k!}{\pi^k(S_n)},\quad k\ge 1,\tag 1.3
$$
so $E_i[T(S_n)]\sim\pi^{-1}(S_n)$ in particular. Thus all the moments of
$T(S_n)/E_i[T(S_n)]$ converge, uniformly over $i$, to the moments of the exponential random variable,
which implies convergence in distribution as well. As for the hit location $X(T(S_n))$,
given $U_n\subseteq S_n$,
$$
\lim_{n\to\infty}\left|P_i\{X(T(S_n))\in U_n\}-\frac{\pi(U_n)}{\pi(S_n)}\right|=0,\tag 1.4
$$
uniformly for $i\in\Bbb N$. Thus,  {\it marginally\/}, $T(S_n)$ and $X(T(S_n))$ behave 
in the limit  as if $X(t)$ is 
a Bernoulli sequence with each trial outcome having distribution $\pi$.
\si

Now suppose that, besides being tight, the chain meets a condition
$$
\delta_0:=\inf_{i,j\in \Bbb N}\sum_{k\in\Bbb N}\,p(i,k)p(j,k)>0.\tag 1.5
$$
For a tight chain, this condition is equivalent to
$$
\rho_0:=\sup_{i,j\in \Bbb N}\,\|p(i,\cdot)-p(j,\cdot)\|_{TV}<1,
$$
which implies that
$$
\|p^n(i,\cdot)-\pi\|_{TV}\le \rho_0^n.
$$
(So, for the relaxation time $\tau$ in [1], we have $\tau=\lceil \ln2/\ln(1/\rho_0)\rceil$.) 

Given a random vector $\bold Y$ with integer components, we denote its probability 
distribution by $d(\bold Y)$. Under the conditions (1.2) and (1.5), we show
that, uniformly for the initial state $i\in \Bbb N$,
$$
\|d((X(T(S_n),T(S_n))-d((\ell_n,t_n))\|_{TV}=O(p(S_n)),\quad p(S):=
\sup_{k\in\Bbb N}p(k,S),\tag 1.6
$$
where $\ell_n$ and $t_n$ are {\it independent\/}, 
$$
\align
P\{t_n=\tau\}=&\,\pi(S_n)(1-\pi(S_n))^{\tau-1},\quad \tau\ge 1,\\
P\{\ell_n=k\}=&\,\frac{\pi(k)}{\pi(S_n)},\quad k\in S_n.
\endalign
$$
More generally, the $k$-long sequence of chronologically ordered locations of first
$k$ hits
of $S_n$, each paired with the time elapsed since a preceding hit (paired with $T(S_n)$ 
in the case of the first hit) is approximated by the $k$-long sequence  of 
{\it independent\/} copies of $(\ell_n,t_n)$, within the total variation distance of order
$O(kp(S_n))$. (Aldous and Brown [2], [3] had used Stein's method to 
show that, for
a stationary, continuous-time, reversible Markov process, the hitting times for a subset $A$ of
states after prolonged excursions outside of $A$ form an approximately Poisson process.)     
\si

The equation (1.6) yields, rather directly, the limiting distributions of the extreme values
for $\{X(t)\}_{1\le t\le N}$. Given $\mu$, let $X^{(\mu)}$ 
denote the $\mu$-th largest among $X(1),\dots, X(N)$. Then
$$
P_i\{X^{(\mu)}\le n\}=P\bigl\{\text{Poisson }(N\pi(S_n))<\mu\bigr\}+O(\mu^2/N+
Np^2(S_n)),\tag 1.7
$$
and we have an extended version of (1.7) for the joint distribution of $X^{(1)},\dots,X^{(\mu)}$.
\bi

Turn now to the application of these results to the random compositions studied
in [19]-[21].
\si

A composition of a positive integer $\nu$ is 
$\bold y=(y_1,\dots,y_\mu)$, $\mu\le \nu$,
such that $y_1,\dots,y_{\mu}$ are positive integers satisfying
$$
\sum_{i=1}^{\mu}y_i=\nu.\tag 1.8
$$
Since, for each $\mu$, there are $\binom{\nu-1}{\mu-1}$ compositions, we have
$2^{\nu-1}$ compositions overall. Assuming that a solution of (1.8) is chosen 
uniformly at random ({\it uar\/}) we have a {\it random\/} composition $\bold Y$ of $\nu$,
its dimension $M$ being random as well. It is known, Andrews [4], 
that 
$$
\bold Y\overset{\Cal D}\to\equiv (Z_1,\dots,Z_{\Cal M-1},\hat Z_{\Cal M}),\tag 1.9
$$
($\overset{\Cal D}\to\equiv$ meaning equality  of distributions), where $Z_1,Z_2,\dots$
are independent Geometrics with success probability $1/2$, 
$$
\Cal M=\min\{m:\,Z_1+\cdots+Z_m\ge \nu\},\tag 1.10
$$
and
$$
\hat Z_{\Cal M}:=\nu-\sum_{j=1}^{\Cal M-1}Z_j.\tag 1.11
$$
Hitczenko and Savage [12] used
this connection to the
well studied  success runs  in a fair coin-tossing process as an efficient tool for 
asymptotic analysis of various characteristics of the random composition.
\si

If a random composition $\bold Y$ is not uniformly distributed on the set
(1.8), one can only hope for {\it asymptotic\/} independence
of most of the parts. Lowering expectations then, one may  search for a Markov chain
that approximates the behavior of $\bold Y$ in question; ergodicity of such a chain would 
mean near independence of parts $Y_{t_1}$ and $Y_{t_2}$ with $|t_1-t_2|$ sufficiently 
large.
\si

Here are two examples of such random compositions. 
A  column-convex-animal (cca) composition of $\nu$ is a collection of lengths of an ordered
sequence of contiguous columns on $\Bbb Z^2$,  whose total sum is $\nu$,
such that every two successive columns
have a common boundary consisting of at least one vertical edge of $\Bbb Z^2$,
Klarner [15], Privman and Forgacs [22], Privman and Svrakic [23], Louchard [19], [20].
A Carlitz (C) composition meets a condition that no two adjacent parts coincide,
Carlitz [8], Knopfmacher and Prodinger [16], Louchard and Prodinger [21], Hitczenko
and Louchard [11]. 

One obtains a certain, {\it nonuniform\/}, distribution on the
set of solutions of (1.8), if a column-convex
animal is chosen uar from among all such creatures.   One obtains another nonuniform distribution,
if a
composition of $\nu$ is chosen uar from among all C-compositions. We call these
objects a random cca composition $\bold Y$, and a random C-composition $\bold Y$, and denote the 
random number of components of $\bold Y$ by $M$. For both schemes, 
Louchard [19], [20] and Louchard and Prodinger [21],
determined  a limiting joint distribution of two successive parts, $Y_t$ and
$Y_{t+1}$, in the case when $t$ and $M-t$ are of order $\nu$, and also the limiting
distribution, $\pi_1$, of the first (last) part $Y_1$ ($Y_M$). These results strongly
suggest, though do not actually prove that, in both cases,
$$
p(i,k):=\lim_{\nu,t\to\infty} \pr\{Y_{t+1}=k\,|\,Y_t=i\},\quad (i,k\in \Bbb N),\tag 1.12
$$
might well be the transition probabilities of a Markov chain, with an initial
distribution $\pi_1$, that closely approximates the {\it whole\/} random composition.
We fully confirm this conjecture, proving an approximational
counterpart of (1.9)-(1.11). The chains turn out to be tight, exponentially
mixing, and this enables us to use our results for asymptotic analysis of
extreme-valued parts of both random compositions. Let $Y^{(\mu)}$ denote the
$\mu$-th largest part of the random composition in question. For the random cca
composition, we show that, for $\mu=o(\ln \nu)$,
$$
Y^{(\mu)}=\frac{\ln\bigl(\mu^{-1}\nu\ln^2 \nu\bigr)}{\ln(1/z^*)}+O_p(1),\tag 1.13
$$
where $z_*=0.31\dots$ is the smallest-modulus root of
$$
4z^3-7z^2+5z-1=0.
$$
For the random C-composition, if $\mu=o(\nu^{1/2})$ then
$$
Y^{(\mu)}=\frac{\ln\bigl(\mu^{-1}\nu\bigr)}{\ln(1/z^*)}+O_p(1),\tag 1.14
$$
where $z_*=0.57\dots$ is the smallest-modulus root of
$$
\sum_{j\ge 1}\frac{z^j}{1+z^j}-1=0.
$$
($O_p(1)$ stands for a random variable bounded in probability.) It follows 
from (1.13) and (1.14) that the number of distinct values among 
$X^{(1)},\dots, X^{(\mu)}$ is likely to be at most $(1+o(1))\ln \mu/\ln(1/z_*)$, for
the corresponding $z_*$, for  $\mu=o(\ln\nu)$ and $\mu=o(\nu^{1/2})$ respectively.
It can be shown that, in fact, the range is asymptotic to $\ln \mu/\ln(1/z_*)$,
in probability. (See Hitczenko and Louchard [11] regarding a limiting
distribution of a ``distinctness''
(range size) of the random C-composition.)
\si

We plan to extend this approach to other constrained compositions, such as 
quite  general Carlitz-type compositions studied by Bender and Canfield [6].
\bi

The rest of the paper is organized as follows. In Section 2 we show that,
for the tight Markov chains $\{X(t)\}$, the hitting time of a rare set $S_n$, i.e. with $\pi(S_n)\to 0$,
scaled by $\pi^{-1}(S_n)$ converges, with all its moments, to the exponentially
distributed random variable of unit mean, while the hit location has, in the limit,
a stationary distribution restricted to $S_n$. And convergence is uniform over all
initial states. In Section 3 we add a second condition
that guarantees exponential mixing, calling such chains tight, exponentially
mixing (t.e.m.) chains. Significantly sharpening the results of Section 2, we demonstrate that
that the hitting time and the hit location are asymptotic, with respect to
the total variation distance, to a pair of {\it independent\/} random variables, one
being geometrically distributed with success probability $\pi(S_n)$, and another
having the restricted stationary distribution. The error term is $O(p(S_n))$,
see (1.6) for definition of $p(\cdot)$. We extend this result to the first $k$
hits of $S_n$, and then state and prove the claims about the limiting distribution of the
$\mu$ largest values among $X(1),\dots, X(N)$, useful for $\mu=o(N^{1/2})$. In
Section 4 we apply these claims to the extreme-valued parts of two random 
compositions of a large $\nu$, the cca composition, and the C-composition. Specifically,
in Section 4.1 we briefly survey the basic known facts about the
compositions. In Section 4.2 we show that each composition is sharply
approximated, in terms of total variation distance, by a related Markov chain,
for as long as the current sum of parts does not exceed $\nu-\ln^2\nu$. In  
Section 4.3, for each composition,  we derive the limiting distributions of the $\mu$ largest
values of a random composition parts, assuming that $\mu=o(\ln\nu)$ for the cca composition,
and $\mu=o(\nu^{1/2})$ for the C-composition. In Appendix we prove an auxiliary
result on large deviations of the number of parts in each of the random compositions.
\bi
\si

{\bf 2. Tight Markov chains.\/} Consider an ergodic Markov chain $X(t)$ on $\Bbb N$
with the stationary distribution $\boldsymbol\pi=\{\pi(j)\}_{j\in \Bbb N}$. Given $S\subset \Bbb N$,
we denote $\pi(S)=\sum_{j\in S}\pi(j)$.  Introduce
$T(S)$ the {\it positive\/} hitting time of $S$, i. e. $T(S)=\min\{t>0\,|\,X(t)\in S\}$, and the
hit location $X(T(S))$.
Our focus is on a {\it rare \/} $S$, i. e. with a small $\pi(S)$. 

Assuming that the chain satisfies a
tightness condition (1), namely
$$
\lim_{n\to\infty}\sup_i\sum_{k>n}p(i,k)=0,\tag 2.1
$$
we will show that, uniformly for an initial state in $\Bbb N$,  (1) $T(S)$ is 
asymptotically exponential, with mean $\pi^{-1}(S)$, and (2) the
distribution of $X(T(S))$ is asymptotic to $\{\pi(s)/\pi(S)\}_{s\in S}$. 
\si

As a first step we prove
the following.

\proclaim{Lemma 2.1} Let a possibly infinite $S_n\neq\emptyset$ be such that 
$\lim_{n\to\infty}\pi(S_n)=0$.
Under the condition (2.1),
$$
E_i[T(S_n)]\sim\frac{1}{\pi(S_n)},\quad n\to\infty,\tag 2.2
$$
uniformly for $i\in \Bbb N$.
\endproclaim
\si

{\bf Note.\/} Consider a simple asymmetric random walk on $\Bbb N$, i. e. the
Markov chain with $p(1,1)=q$,
$p(1,2)=p$, and $p(i,i-1)=q$, $p(i,i+1)=p$ for $i\ge 2$. For $p<q$ this chain is ergodic,
with the stationary distribution $\pi(j)=(1-p/q)(p/q)^{j-1}$, but it is clearly not
tight. For $i=1$, $T(\{n+1\})=T(\{n+1,n+2,\dots\})$, but
$\pi(\{n+1\})\not\sim \pi(\{n+1,n+2,\dots\})$. So (2.2) cannot hold for all $S_n$
with $\pi(S_n)\to 0$. In fact, the expected common hitting time for these two sets
is not asymptotic to the reciprocal of either of these stationary probabilities. 

\bi
{\bf Proof of Lemma 2.1.\/}  By tightness condition (2.1), there exists $K$ such that
$$
\sum_{j\le K}p(i,j)\ge 1/2,\quad \forall i\ge 1.
$$
Then, for $t\ge 1$,
$$
P_i\{T([K])> t\} \le \frac{1}{2^t}\Longrightarrow E_i[T([K])]\le 2.
$$
Now, one (possibly not the shortest) way of hitting $S_n$, starting at $i$,
is to hit
the set $[K]$ and from there to hit $S_n$. By the strong Markov
property, conditionally on $X(T([K]))=j,\,(j\in [K])$, the residual travel time
$\hat T(S_n)$ till hitting $S_n$ is distributed as $T(S_n)$ under $P_j$. So 
$$
E[\hat T(S_n)|X(T([K]))=j] =E_j[T(S_n)],\quad j\in [K].
$$
Then, introducing $\ell\in [K]$ such that
$$
E_{\ell}[T(S_n)]=\max_{j\in [K]}E_j[T(S_n)],
$$
we have: 
$$
\aligned
E_i[T(S_n)]\le& E_i[T[K]] +\sum_{j\in [K]}P_i\{X(T([K]))=j\}E_j[T(S_n)]\\
\le& 2 + E_{\ell}[T(S_n)];
\endaligned\tag 2.3
$$
in particular, $\sup_iE_i[T(S_n)]<\infty$.
\si

By Markov property, 
$$
E_j[T(S_n)]=1+\sum_{k\in S_n^c}p(j,k)E_k[T(S_n)],\quad j\in \Bbb N.\tag 2.4
$$
Multiplying both sides of (2.4) by $\pi(j)$ and summing for $j\in \Bbb N$, we get
$$
\align
\sum_{j\in \Bbb N}\pi(j)E_j[T(S_n)]=&\,1+
\sum_{k\in S_n^c}E_k[T(S_n)]\sum_{j\in \Bbb N}\pi(j)p(j,k)\\
=&\,1+\sum_{k\in S_n^c}\pi(k)E_k[T(S_n)],
\endalign
$$
as $\pi(\cdot)$ is stationary. So, as both series converge,
$$
\sum_{k\in S_n}\pi(k)E_k[T(S_n)]=1. \tag 2.5
$$
(We note that (2.5) is a special case of a well-known result, due to Kac [13],  with inevitably harder
proof,  for a general discrete-time stationary process; see also Breiman [7],
Section 6.9.) Then, by (2.3),
$$
E_{\ell}[T(S_n)]+1\ge \frac{1}{\pi(S_n)}\Longrightarrow E_{\ell}[T(S_n)]\gtrsim\frac{1}{\pi(S_n)}. \tag 2.6
$$

Now, given a state $k$, we have
$$
E_{\ell}[T(S_n)]\le E_{\ell}[T(\{k\})] +E_k[T(S_n)], \tag 2.7
$$
$T(\{k\})$ being the hitting time for the singleton $\{k\}$.
Combining (2.6) and (2.7), we obtain: for every {\it fixed\/} $k$,
$$
E_k[T(S_n)] \gtrsim \frac{1}{\pi(S_n)}.\tag 2.8
$$
Picking  arbitrary $L$, by (2.4), we have: for $n\ge n(L)$,
$$
E_j[T(S_n)]\ge 1 +\sum_{k\le L}p(j,k)E_k[T(S_n)], \quad j\in \Bbb N.
$$
Therefore, by (2.8),
$$
\liminf_{n\to\infty} \left(\inf_{j\in \Bbb N} E_j[T(S_n)]\right)\pi(S_n)
\ge \liminf_{n\to\infty}\inf_{j\in \Bbb N}\sum_{k\le L}p(j,k),
$$
where, by (2.1), the RHS approaches $1$ as $L\uparrow\infty$. So
$$
E_j[T(S_n)]\gtrsim\frac{1}{\pi(S_n)},\tag 2.9
$$
{\it uniformly\/} for $j\in \Bbb N$.
\si

It remains to show that 
$$
E_j[T(S_n)]\lesssim \frac{1}{\pi(S_n)},
$$
uniformly for $j\in \Bbb N$.  Using (2.4)-(2.5), we obtain then
$$
\sum_{j\in S_n}\frac{\pi(j)}{\sum\limits_{i\in S_n}\pi(i)}\left(1+\sum_{k\in S_n^c}p(j,k)
E_k[T(S_n)]\right)=
\frac{1}{\pi(S_n)}.\tag 2.10
$$
Suppose that there exists a subsequence $n_m\to\infty$ and $\delta>0$, such that
$$
\lim_{n\in \{n_m\}}E_{\ell}[T(S_n)]\pi(S_n)\ge 1+\delta.
$$
Then, by (2.7),
$$
\lim_{n\in \{n_m\}}E_{k}[T(S_n)]\pi(S_n)\ge 1+\delta,
$$
for every fixed $k$. Picking $M>0$ and dropping the summands for $k>M$ in (2.10), we get
then: for $n=n_m$ large enough,
$$
(1+\delta/2)\sum_{j\in S_n}\frac{\pi(j)}{\sum\limits_{i\in S_n}\pi(i)}\left(\sum_{k\le M}
p(j,k)\right)\le 1.
$$
This is impossible if $M$ is chosen so large that
$$
\inf_j\sum_{k\le M}p(j,k)\ge \frac{1}{1+\delta/3}.
$$
Therefore
$$
E_{\ell}[T(S_n)]\lesssim \frac{1}{\pi(S_n)},
$$
and so, invoking (2.3),
$$
E_k[T(S_n)]\lesssim \frac{1}{\pi(S_n)}, \tag 2.11
$$
uniformly for $k\in \Bbb N$.
\si
 
Combining (2.9) and (2.11), we complete the proof of Lemma 2.1.\qed
\bi

The fact that $E_i[T(S_n)]\to\infty$  already implies, via
Keilson's theorem [14], that, for each {\it fixed\/} initial state $i$, $T(S_n)
/E_i[T(S_n)]$ is, in the limit, exponentially distributed, with parameter $1$. The tightness
condition allowed us to estimate the scaling parameters $E_i[T(S_n)]$ asymptotically,
uniformly for $i\in \Bbb N$. Interestingly, this uniformity can be used for a simple
alternative proof of asymptotic exponentiality of $T(S_n)/E_i[T(S_n)]$.
\proclaim{Lemma 2.2} Under the condition (2.1), for each fixed $k\ge 1$,
$$
E_i[T^k(S_n)]\sim k!/\pi^k(S_n),\tag 2.12
$$
uniformly for $i\in \Bbb N$. Consequently, uniformly for $i\in \Bbb N$,
$$
P_i\{T(S_n)\pi(S_n)>x\}\to e^{-x},\quad\forall\,x\ge 0.\tag 2.13
$$
\endproclaim
\si
{\bf Proof of Lemma 2.2.\/} Introduce the moment generating functions
$$
\phi_i(u)=\sum_{r\ge 0}\frac{u^r}{r!}E_i[T^r(S_n)],\quad i\in \Bbb N.
$$
As formal power series, these functions satisfy
$$
\phi_i(u)=e^u\left(\sum_{j\in S_n}p(i,j)+\sum_{j\in S_n^c}p(i,j)\phi_j(u)\right),
\quad i\in \Bbb N.\tag 2.14
$$ 
Differentiating both sides of (2.14) $k$ times at $u=0$ we get
$$
\aligned
E_i[T^k(S_n)]=\,& b(i,k)
+\sum_{j\in S_n^c}p(i,j)E_j[T^k(S_n)],\quad i\in \Bbb N,\\
b(i,k):=&1+\sum_{r=1}^{k-1}\binom{k}{r}\sum_{j\in S_n^c}p(i,j)E_j[T^r(S_n)].
\endaligned\tag 2.15
$$
For $k=1$ we get (2.4). Let $k\ge 2$, and suppose that, for $r<k$,
$$
E_i[T^{r}(S_n)]=(1+o(1))\frac{r!}{\pi^{r-1}(S_n)}E_i[T(S_n)],\tag 2.16
$$
uniformly for $i\in \Bbb N$. (This is obviously true for $k=2$.) Then
$$
\align
b(i,k)=\,&1 +(1+o(1))\sum_{r=1}^{k-1}\frac{(k)_r}{\pi^{r-1}(S_n)}
\sum_{j\in S_n^c}p(i,j)E_j[T(S_n)]\\
&(\text{using }(2.4))\\
=\,&1 +(1+o(1))\sum_{r=1}^{k-1}\frac{(k)_r}{\pi^{r-1}(S_n)}\bigl(E_i[T(S_n)]-1)\\
=\,&(1+o(1))\frac{(k)_{k-1}}{\pi^{k-1}(S_n)}=(1+o(1))\frac{k!}{\pi^{k-1}(S_n)},\tag 2.17
\endalign
$$
uniformly for $i\in \Bbb N$. Using  (2.17), we rewrite (2.15) as
$$
E_i[T^k(S_n)]=(1+o(1))\,\frac{k!}{\pi^{k-1}(S_n)}+\sum_{j\in S_n^c}p(i,j)E_j[T^k(S_n)], 
$$
uniformly for  $i\in \Bbb N$. Now if define $x^0(i,k)=b(i,k)$, and, for $t\ge 0$,
$$
x^{t+1}(i,k)=b(i,k)+\sum_{j\in S_n^c}p(i,j)x^t(j,k),\quad i\in \Bbb N,
$$
then $x^t(i,k)\uparrow E_i[T^k(S_n)]$, $i\in \Bbb N$. In particular, for $k=1$,
we have $b(i,1)=1$, and $x^t(i,1)\uparrow E_i[T(S_n)]$. Using this observation and
(2.17), and $E_i[T(S_n)]\sim 1/\pi(s_n)$
uniformly for $i\in\Bbb N$, we conclude:
$$
E_i[T^k(S_n)]=(1+o(1))\frac{k!}{\pi^{k-1}(S_n)}\,E_i[T(S_n)],
$$
uniformly for $i\in \Bbb N$. Thus (2.16) holds for all $r\ge 1$, and so
$$
E_i[T^k(S_n)]=(1+o(1))\frac{k!}{\pi^k(S_n)},\quad k\ge 1,
$$
uniformly for $i\in \Bbb N$.
\si

Since $\limsup k^{-1} (k!)^{1/k}/k<\infty$, the exponential distribution is the only one with the moments $k!$,
(Durrett [10]). The proof of Lemma 2.2 is complete.\qed
\bi

Turn now to $H_n:=X(T(S_n))$, $H$ reminding us that $X(T(S_n))$ is the {\it hit\/} location.
\proclaim{Lemma 2.3} Let $U_n\subseteq S_n$. Uniformly for $i\in \Bbb N$,
$$
\lim_{n\to\infty}\left|P_i\{H_n\in U_n\}-\frac{\pi(U_n)}{\pi(S_n)}\right|=0.\tag 2.18
$$
\endproclaim
\si 
{\bf Proof of Lemma 2.3.\/} By Markov property,
$$
P_i\{H_n\in U_n\}=p(i,U_n)+\sum_{j\in S_n^c}p(i,j)P_j\{H_n\in  U_n\},\quad
i\in \Bbb N,\tag 2.19
$$ 
where we use the notation $p(i,A)=\sum_{k\in A}p(i,k)$, $A\subseteq \Bbb N$. 

(a) Assuming only  that $\{p(i,k)\}$ is ergodic,
let us show that, for all {\it fixed\/} $i,j\in \Bbb N$, 
$$
\lim_{n\to\infty} \bigl|P_i\{H_n\in U_n\}-P_j\{H_n\in U_n\}\bigr|=0.
\tag 2.20
$$
By Cantor diagonalization device, any
subsequence $\{n_m\}$ of $1,2,\dots$ contains a further subsequence 
$\{n_{m_{\ell}}\}$
such that for $n\to\infty$ along this subsequence,  there exists
$$
f_i=\lim_{n\to\infty}P_i\{H_n\in U_n\}\in [0,1],\quad i\in \Bbb N
$$
The limits $f{_i}$ may well depend on $\{n_m\}$, of course. 
Letting $n=n_{m_{\ell}}\to\infty$ in (2.19), we obtain: 
$$
f_i=\sum_{j\in \Bbb N}p(i,j) f_{j},\quad i\in \Bbb N.\tag 2.21
$$
Since the matrix $\{p(i,j)\}$ is ergodic, $f_{i}$ does not depend
on $i$, (Durrett [10], Exer. 3.9). So (2.20) follows.

(b) By tightness, given $\varepsilon\in (0,1)$, there exists
$J=J(\varepsilon)$ such that
$$
\sum_{j\le J}p(i,j)\ge 1-\varepsilon,\quad \forall\, i\in \Bbb N.
$$
For $n\ge n(J)$, $[J]\subseteq S_n^c$.  So, by (2.19), 
$$
\inf_{i\in \Bbb N}P_i\{H_n\in U_n\}\ge (1-\varepsilon)\min_{i\le J}
P_i\{H_n\in U_n\},
$$
and
$$
\sup_{i\in \Bbb N}P_i\{H_n\in U_n\}\le \varepsilon +\max_{i\le J}
P_i\{H_n\in U_n\}.
$$
So
$$
\multline
\limsup_n\left[\sup_{i\in \Bbb N}P_i\{H_n\in U_n\}-\inf_{i\in \Bbb N}P_i\{H_n\in U_n\}
\right]\\
\le 2\varepsilon+\lim_{n\to\infty}\bigl[\max_{i\le J}P_i\{H_n\in U_n\}-\min_{i\le J}
P_i\{H_n\in U_n\}\bigr]
=2\varepsilon.
\endmultline
$$
Thus
$$
\lim_{n\to\infty}\left[\sup_{i\in \Bbb N}P_i\{H_n\in U_n\}-\inf_{i\in \Bbb N}P_i\{H_n\in U_n\}
\right]=0.\tag 2.22
$$

(c) Multiplying both sides of (2.19) by $\pi(i)$,  summing for $i\in \Bbb N$, and using
stationarity of $\pi(\cdot)$, we obtain
$$
\align
\sum_{i\in \Bbb N}\pi(i)P_i\{H_n\in U_n\}=&\,\sum_{i\in \Bbb N}\pi(i)p(i,U_n)
+\sum_{j\in S_n^c}P_j\{H_n\in U_n\}\sum_{i\in\Bbb N}\pi(i)p(i,j)\\
=&\,\pi(U_n)+\sum_{j\in S_n^c}\pi(j)P_j\{H_n\in U_n\},
\endalign
$$
so that
$$
\sum_{i\in S_n}\pi(i)P_i\{H_n\in U_n\}=
\pi(U_n).\tag 2.23
$$
Consequently
$$
\inf_{i\in S_n}P_i\{H_n\in U_n\}\le q(n,U_n)\le \sup_{i\in S_n}
P_i\{H_n\in U_n\},
$$
where
$$
q(n,U_n):=\frac{\pi(U_n)}{\pi(S_n)}.
$$
Combining (2.22) and the double inequality we conclude that
$$
\lim_{n\to\infty}\left|P_i\{H_n\in U_n\}-\frac{\pi(U_n)}{\pi(S_n)}\right|=0,
$$
uniformly for $i\in \Bbb N$. The proof of Lemma 2.3 is complete. \qed
\bi

Thus, considered separately, $T(S_n)$ and $X(T(S_n))$ asymptotically behave
as if $X(t)$ is a Bernoulli sequence with each trial outcome
having distribution $\pi$. Of course, the Bernoulli
sequence possesses finer properties; in particular, $T(S_n)$ and $X(T(S_n))$
are independent of each other. We are about to impose an additional
condition on $\{p(i,k)\}$.  It will be used to to establish a limit distribution of the
vector $\bigl(T(S_n),X(T(S_n))\bigr)$, together with a convergence rate in terms
of the $\|\cdot\|_{TV}$ distance.
In particular,  under the two conditions,  $T(S_n))$ and 
$X(T(S_n))$ turn out to be {\it asymptotically\/} independent.
\bi
{\bf 3. Tight, exponentially mixing Markov  chains.\/} The  extra condition (2) 
is:
$$
\rho:=\sup_{i,j\in\Bbb N}\sum_{k\in \Bbb N}|p(i,k)-p(j,k)|<2. \tag 3.1
$$
(Of course, $\rho\le 2$ always.) Then (Durrett [10], Exer. 5.11), 
$$
\sum_{k\in \Bbb N}|p^n(i,k)-p^n(j,k)|\le 2 (\rho/2)^n,
$$
where $p^n(\cdot,\cdot)$ are the $n$-step transition probabilities. Consequently,
multiplying by $\pi(j)$ and summing over $j\in \Bbb N$,
$$
\sum_{k\in \Bbb N}|p^n(i,k)-\pi(k)|\le 2 (\rho/2)^n. 
$$
Equivalently, denoting  $\bold e=
(\{1\}_{i\in\Bbb N})^T$,
$$
\|\bigl(P^n - \bold e\boldsymbol\pi\bigr)^T\|_{L_1(\Bbb N)}
=\|P^n - \bold e\boldsymbol\pi\|_{L_{\infty}(\Bbb N)}\le 2 (\rho/2)^n.
\tag 3.2
$$
We call the  chains meeting (3.2) exponentially mixing, and use abbreviation
t. e. m. chains for tight, exponentially mixing Markov chains.
Now 
$$
\frac{1}{2}\sum_{k\in \Bbb N}|p(i,k)-p(j,k)|=\min_{(X,Y)}P\{X\neq Y\},
$$
where minimum is over all random vectors $(X,Y)$ such that $P\{X=k\}=p(i,k)$,
$P\{Y=k\}=p(j,k)$, $k\in \Bbb N$, see Durrett [10]. Therefore, selecting independent $X$ and $Y$,
$$
\frac{1}{2}\sum_{k\in \Bbb N}|p(i,k)-p(j,k)|\le 1-P\{X=Y\}=1-
\sum_{k\in \Bbb N}p(i,k)p(j,k).
$$
Hence the condition (3.1) is met if
$$
\delta_0:=\inf_{i,j\in \Bbb N}\sum_{k\in \Bbb N}p(i,k)p(j,k)>0,\tag 3.3
$$
in which case $\rho/2\le 1-\delta_0$. In fact, for the tight chains the converse is true: 
(3.1) implies (3.3). Suppose not. Then there exists $\{(i_r,j_r)\}_{r\ge 1}$ such that
$$
\lim_{r\to\infty}\sum_{k\in \Bbb N}p(i_r,k)p(j_r,k)=0.\tag 3.4
$$
By the tightness condition, we may assume that $p(i_r,\cdot)$ and $p(j_r,\cdot)$
converge, weakly, to some probability distributions, $p_1$ and $p_2$ respectively,
that is
$$
p(i_r,k)\to p_1(k),\quad p(j_r,k)\to p_2(k),\quad k\in \Bbb N.
$$
Combining this with (3.1) and (3.4), we obtain
$$
\sum_{k\in \Bbb N}|p_1(k)-p_2(k)|<2,\quad \sum_{k\in \Bbb N}p_1(k)p_2(k)=0.
$$
This is impossible, since the second condition implies that
$$
|p_1(k)-p_2(k)|=p_1(k)+p_2(k).
$$
\bi

\proclaim{Theorem 3.1} Let $\partial_{n,i}$ denote the joint distribution
of $X(T(S_n))$ and $T(S_n)$ for  an initial state $i\in\Bbb N$. Let $\partial_n$
denote the product probability measure on $S_n\times \Bbb N$, such that
$$
\partial_n(A\times B) =\frac{\pi(A)}{\pi(S_n)}\cdot\sum_{\tau\in B}\pi(S_n)(1-\pi(S_n))
^{\tau-1},\quad A\subseteq S_n,\,B\subseteq \Bbb N.
$$
Under the conditions (1) and (2), uniformly for $i\in \Bbb N$, 
$$
\|\partial_{n,i}-\partial_n\|_{TV}=O(p(S_n)),\tag 3.5
$$
where 
$$
p(S_n):=\sup\limits_{i\in \Bbb N}p(i,S_n).\tag 3.6
$$
\endproclaim
\si
{\bf Proof of Theorem 3.1.\/} Introduce
$$
\e_n=\sup_{i\in S_n^c}p(i,S_n),\quad p(i,A):=\sum_{k\in A}p(i,k);
$$
by the tightness and $\lim\pi(S_n)=0$, we have $\lim\e_n=0$. Let $P_n=\{p(i,k)\}_{i,k\in
S_n^c}$. 
\si

As a first step let us prove the following claim.
\proclaim{Lemma 3.2} For $n$ large enough, $P_n$ has an eigenvalue $\lambda_n\in [1-\e_n,1)$ 
and
a corresponding eigenvector $\bold f_n=\left(\{f_n(i)\}_{i\in S_n^c}\right)^T$, such that
$$
1\le f_n(i)\le\frac{1}{1-(6/\delta_0)\e_n},\tag 3.7
$$
with $\delta_0$ coming from the condition (3.3).
\endproclaim
\si
{\bf Proof of Lemma 3.2.\/} Given $m>0$, introduce
$S_{n,m}=S_n\cup \{m+1,m+2,\dots\}$, so that $S_{n,m}^c=S_n^c\cap [m]$,
which is a finite set. Denote $P_{n,m}=\{p(i,k)\}_{i,k\in S_{n,m}^c}$.
By the conditions (2.1) and (3.3), there exist $n_0$ and $m_0$ such that
$$
\align
&\inf_{n\ge n_0\atop m\ge m_0}\min_{i\in S_{n,m}^c}p(i,S_{n,m}^c)
>1-
\delta_0/3,\tag 3.8\\
&\inf_{n\ge n_0\atop m\ge m_0}\min_{i,j\in S_{n,m}^c}\sum_{k\in S_{n,m}^c}
p(i,k)p(j,k)\ge \delta_0/2. \tag 3.9
\endalign
$$
Let $n\ge n_0$, $m\ge m_0$. Call $\emptyset\neq A\subseteq S_{n,m}^c$ closed (in $S_{n,m}^c$), if 
$p(i, S_{n,m}^c\setminus A)=0$ for each $i\in A$. Call
a closed set minimal if it does not contain a closed subset. The
condition (3.9) clearly ensures that there exists exactly one  minimal
closed subset $A$, which may be the whole set $S_{n,m}^c$. A submatrix
$P_A:=\{p(i,k)\}_{i,k\in A}$ is irreducible; so it has a positive eigenvalue 
$\lambda(A)$ with a positive eigenvector $\bold f_A$, and 
the absolute values of the remaining eigenvalues of $P_A$ do not exceed 
$\lambda(A)$. In fact, those absolute values are strictly less than $\lambda(A)$.
Otherwise, by Frobenius theorem, there exists a partition $A=\uplus_{r=1}^h A_r$, $h>1$,
such that, for $r=1,\dots,h$, $\{k\in A\,|\,\exists\,i\in A_r,\,p(i,k)>0\}=A_{r+1}$, 
($h+1:=1$). So, for $i\in A_1$, $j\in A_2$,
$$
\sum_{k\in L}p(i,k)p(j,k)=\sum_{k\in A}p(i,k)p(j,k)=\sum_{k\in A_2\cap A_3}p(i,k)p(j,k)=0.
$$
And this contradicts (3.9).  Furthermore, by (3.8),
$$
\lambda(A)\ge \min_{i\in A}p(i,A)=\min_{i\in A}p(i,S_{n,m}^c)>1-\delta_0/3,
$$
while, denoting $A^\prime=S_{n,m}^c\setminus A$ and using (3.9),
$$
\align
\max_{i\in A^\prime}p(i,A^\prime)\le\, & 
1-\min_{i\in A^\prime}p(i,A)\\
\le\, &1-\min_{i\in A^\prime,\, j\in A}\sum_{k\in A}p(i,k)p(j,k)\\
\le\,& 1-\min_{i\in A^\prime,\,j\in S_{n,m}^c}\sum_{k\in S_{n,m}^c}p(i,k)p(j,k)
\le1-\delta_0/2.
\endalign
$$
Therefore $\lambda(A)$ is strictly larger than $\lambda(A^\prime)$, the largest
eigenvalue of  $P_{A^\prime}$. Denoting $P_{A^\prime, A}=\{p(i,k)\}_{i\in A^\prime,
k\in A}$, let $\bold f_{A^\prime}$ be a solution of
$$
P_{A^\prime,A}\bold f_A +P_{A^\prime}\bold f_{A^\prime}=\lambda(A)
\bold f_{A^\prime}.
$$
Since $P_{A^\prime,A}\bold f^A>\bold 0$ and $\lambda(A)>\lambda(A^\prime)$,
$\bold f_{A^\prime}$ exists uniquely and is positive. The combined vector $\bold f_{n,m}=
(\bold f_A,\bold f_{A^\prime})$ is a unique, positive, eigenvector of $P_{n,m}$ for 
its
largest eigenvalue 
$$
\lambda_{n,m}=\lambda(A)\ge 1-\delta_0/3.
$$
\si

Let us bound $\max_i\bold f_{n,m}(i)/\min_i\bold f_{n,m}(i)$. Introduce 
$p_i(\tau)$, ($i\in S_{n,m}^c,\,\tau\ge 0$),
the probability that, starting at state $i$, the Markov process $X(t)$ stays in 
$S_{n,m}^c$ for all $t\le \tau$.
The sequence $\bold p(\tau)=\{p_i(\tau)\}_{i\in S_{n,m}^c}$, satisfies a 
recurrence
$$
\bold p(\tau+1)=P_{S_{n,m}^c}\bold p(\tau),\quad \bold p(0)=(1,\dots,1)^T.
$$
Moreover, there exists $C_{n,m}>0$ such that
$$
\bold p(\tau)\sim C_{n,m}\,\lambda^{\tau}_{n,m}\, \bold f_{n,m},\quad \tau\to\infty.
\tag 3.10
$$
To exploit this connection, let us first use a coupling device to derive a recurrence for
the differences $p_i(\tau)-p_j(\tau)$, $i\neq j$. Consider two independent processes, $X(t)$ and $Y(t)$,
starting at $i$ and $j$ in $S_{n,m}^c$. Introduce the events 
$U(t)=\{X(t)\in S_{n,m}^c\}$ and
$V(t)=\{Y(t)\in S_{n,m}^c\}$, and let $\bold 1(W)$ denote the indicator of an event 
$W$.
By the Markov property,
$$
\aligned
&p_i(\tau+1)-p_j(\tau+1)
=\,E_{(i,j)}\left[\prod_{t\le \tau+1}\bold 1(U(t))-\prod_{t\le \tau+1}\bold 1(V(t))
\right]\\
=&\sum_{k_1,k_2\in S_{n,m}^c\atop k_1\neq k_2}p(i,k_1)p(i,k_2)
E_{(k_1,k_2)}\left[\prod_{t\le \tau}\bold 1(U(t))-\prod_{t\le \tau}\bold 1(V(t))
\right]\\
&+\sum_{k_1\in S_{n,m}^c,\, k_2\in S_{n,m}}
p(i,k_1)p(i,k_2)E_{k_1}\left[\prod_{t\le \tau}\bold 1(U(t))\right]\\
&-\sum_{k_1\in S_{n,m},\, k_2\in S_{n,m}^c}p(i,k_1)p(i,k_2)E_{k_2}
\left[\prod_{t\le \tau}\bold 1(V(t))\right]\\
=&\sum_{k_1,k_2\in S_{n,m}^c\atop k_1\neq k_2}p(i,k_1)p(i,k_2)[p_{k_1}(\tau)-p_{k_2}(\tau)]\\
&+p(j,S_{n,m})\!\!\!\!\sum_{k_1\in S_{n,m}^c}\!\!\!\!p(i,k_1)p_{k_1}(\tau)\,\,
-\,\,p(i,S_{n,m})\!\!\!\sum_{ k_2\in S_{n,m}^c}\!\!\!p(j,k_2)p_{k_2}(\tau).
\endaligned
$$
Letting $\tau\uparrow\infty$ and using (3.10) we obtain 
$$
\multline
\lambda_{n,m}(f_{n,m}(i)-f_{n,m}(j))=\sum_{k_1,k_2\in S_{n,m}^c\atop k_1\neq k_2}
\!\!\!\!p(i,k_1)p(j,k_2)(f_{n,m}(k_1)-f_{n,m}(k_2))\\
+p(j,S_{n,m})\!\!\!\sum_{k_1\in S_{n,m}^c}\!\!\!p(i,k_1)f_{n,m}(k_1)\,\,
-\,\,p(i,S_{n,m})\!\!\!\sum_{ k_2\in S_{n,m}^c}\!\!\!p(j,k_2)f_{n,m}(k_2).
\endmultline\tag 3.11
$$
Let $f_{n,m}(i_1)=\max_{i\in S_{n,m}^c}f_{n,m}(i)$, $f_{n,m}(i_2)=\min_{i\in S_{n,m}^c}
f_{n,m}(i)$. Then it follows from (3.11) that
$$
\bigl(f_{n,m}(i_1)-f_{n,m}(i_2)\bigr)\left(\lambda_{n,m}-\sum_{k_1,k_2\in S_{n,m}^c
\atop k_1\neq k_2}\!\!\!\!p(i_1,k_1)p(i_2,k_2)\right)\le
p(i_2,S_{n,m})\,f_{n,m}(i_1).
$$
Here, see (3.3),
$$
\sum_{k_1,k_2\in S_{n,m}^c\atop k_1\neq k_2}p(i_1,k_1)p(i_2,k_2)\le\,
1-\sum_{k\in S_{n,m}^c}p(i_1,k)p(j,k)\le 1-\delta_0/2,
$$
and
$$
p(i_2,S_{n,m})\le \e_{n,m}:=\max_{i\in S_n^c}\,p(i,S_{n,m}).
$$
As $\lambda_{n,m}\ge 1-\delta_0/3$, we obtain
$$
f_{n,m}(i_1)\left(1-\frac{6}{\delta_0}\e_{n,m}\right)\le f_{n,m}(i_2),\quad 
\forall\, n\ge n_0,\,\,\forall m\ge m_0.\tag 3.12
$$
Now 
$$
0\le \e_{n,m}-\e_n\le\sup_i\sum_{k>m}p(i,k),
$$
so that $\lim_{m\to\infty}\e_{n,m}=\e_n$, uniformly for $n$, and $\lim_{n\to\infty}\e_n=0$.
So there exist
$n_1>n_0$, $m_1\ge m_0$ such that $\e_{n,m}\le \delta_0/7$ for $n\ge n_1$, $m\ge m_1$.
For those $n,m$, the relation (3.12), with $f_{n,m}(i_2)=\min_if_{n,m}(i)=1$, yields
$$
1\le f_{n,m}(i)\le \frac{1}{1-(6/\delta_0)\e_{n,m}},\quad i\in S_{n,m}^c. \tag 3.13
$$
A standard argument shows then existence of a subsequence $m_s\uparrow\infty$ such that (1) 
for each $i\in S_n^c$, there exists $f_n(i)=\lim_{m_s\to\infty}f_{n,m_s}(i)$, which
necessarily satisfies
$$
1\le f_n(i)\le\frac{1}{1-(6/\delta_0)\e_n},\quad i\in S_n^c,
$$
and (2) there exists $\lambda_n=\lim_{m_s\to\infty}\lambda_{n,m_s}\in [1-\e_n,1)$.
Clearly then $\bold f_n:=\bigl(\{f_n(i)\}_{i\in S_n^c}\bigr)^T\in L_{\infty}(S_n^c)$ is an 
eigenvector of $P_n$, $\lambda_n$  being a corresponding eigenvalue. 
The proof of Lemma 3.2 is complete.\qed 
\bi

Let $F_n$ be a diagonal $S_n^c\times S_n^c$ matrix with $F_n(i,i)=f_n(i)$, $i\in S_n^c$.
Define a $S_n^c\times S_n^c$ matrix
$$
Q_n=\lambda_n^{-1}F_n^{-1}P_nF_n=\lambda_n^{-1}\{(f_n(i))^{-1}p(i,k)f_n(k)\}_{i,k\in S_n^c}.
$$
Let $\bold e_n=\bigl(\{1\}_{i\in S_n^c}\bigr)^T$. Since $F_n\bold e_n=\bold f_n$, we have
$$
Q_n\bold e_n=\lambda_n^{-1}F_n^{-1}P_n\,\bold f_n=F_n^{-1}\bold f_n=\bold e_n,
$$
so that $Q_n$ is stochastic. From tightness of $P=\{p(i,k)\}_{i,k\in \Bbb N}$,
and (3.9) it follows that, for each fixed $n$, and even uniformly over
$n$, $Q_n$ is tight as well, i. e.
$$
\delta(K):=\sup_{n,\,i\in S_n^c}\sum_{k\in S_n^c: k>K}Q_n(i,k)\to 0,\quad K\uparrow
\infty.
$$
Now from
$$
(Q_n)^{\nu}(i,k)=\sum_{j\in S_n^c}Q_n(i,j)(Q_n)^{\nu-1}(j,k),
$$
by induction on $\nu$ it follows that
$$
\sup_{\nu,\,i\in S_n^c}\sum_{k\in S_n^c:\, k>K}(Q_n)^{\nu}(i,k)\le \delta(K).
$$
Hence, given $n$, the rows of all matrices $(Q_n)^{\nu}$ form a tight set of
probability distributions. 
Therefore there exists $\nu_s\to\infty$ and a family of probability
distributions $\pi_n(i,\cdot)$ on $S_n^c$, ($i\in S_n^c$), such that, for $i,k\in S_n^c$,
$$
(Q_n)^{\nu}(i,k)\to \pi_n(i,k),\quad \nu\to\infty.
$$ 
In addition, by (3.7) and (3.9), for $i,j\in S_n^c$,
$$
\align
\sum_{k\in S_n^c}Q_n(i,k)Q_n(j,k)=&\,\lambda_n^{-2}\sum_{k\in S_n^c}
\frac{(f^n_k)^2}{f^n_if^n_j}p(i,k)p(j,k)\\
\ge&\,\bigl(1-(6/\delta_0)\e_n\bigr)^2\delta_0/2\ge \delta_0/3,
\endalign
$$
for $n$ large enough. Therefore, cf. (3.2),
$$
\sum_{k\in S_n^c}|(Q_n)^{\nu}(i,k)-(Q_n)^{\nu}(j,k)|\le 2(1-\delta_0/3)^{\nu}\to 0,
\quad\nu\to\infty.
\tag 3.14
$$
Letting $\nu\to\infty$ along $\{\nu_s\}$ in (3.14) we obtain that the family
$\{\pi_n(i,\cdot)\}_{i\in S_n^c}$ consists of a single probability 
distribution $\pi_n(\cdot)$ on $S_n^c$. Thus, for any distribution $q(\cdot)$ on 
$S_n^c$, $q(Q_n)^{\nu_s}\to \pi_n$. Applying this to $q=\pi_n$, and then to 
$q=\pi_nQ_n$,
$$
\pi_nQ_n=\lim_{\nu_s\to\infty}\pi_n(Q_n)^{\nu_s}Q_n=\lim_{\nu_s\to\infty}
(\pi_nQ_n)(Q_n)^{\nu_s}=\pi_n,
$$
that is $\pi_n$ is a stationary distribution of $Q_n$. Using stationarity of
$\pi_n$ and (3.14) we obtain 
$$
\bigl\|\bigl[(Q_n)^{\nu}-\bold e_n\boldsymbol\pi_n\bigr]^T\bigr\|_{L_1(S_n^c)}=
\bigl\|(Q_n)^{\nu}-\bold e_n\boldsymbol\pi_n\bigr\|_{L_{\infty}(S_n^c)}\le 
2(1-\delta_0/3)^{\nu},\tag 3.15
$$
cf. (3.2). Since $(Q_n)^{\nu}=\lambda_n^{-\nu}F_n^{-1}P_n^{\nu}F_n$, combination
of (3.3) and (3.2) implies that
$$
(P_n)^{\nu}=\,\lambda_n^{\nu}\,\bold f_n \boldsymbol{\sigma}_n+R_{n,\nu},\quad
\boldsymbol\sigma_n:=\{\pi_n(i)/f_n(i)\}_{i\in S_n^c}, \tag 3.16
$$
where
$$
\bigl\|R_{n,\nu}^T\bigr\|_{L_1(S_n^c)}=\|R_{n,\nu}\|_{L_{\infty}(S_n^c)}\le
2(1-\delta_0/3)^{\nu}.\tag 3.17
$$
\bi

The estimates (3.16)-(3.17) enable us to determine the limiting joint
distribution of $T(S_n)$ and $X(T(S_n))$. Given $A\subseteq S_n$, and
$z$ with $|z|\le 1$, define 
$$
\psi_i(z)=E_i\bigl[z^{T(S_n)}\bold 1(X(T(S_n))\in A)\bigr],\quad i\in S_n^c,
$$
and $\boldsymbol\psi(z)=\bigl[\{\psi_i(z)\}_{i\in S_n^c}\bigr]^T$. Using Markov 
property, we have
$$
\psi_i(z)=zp(i,A)+z\sum_{k\in S_n^c}p(i,k)\psi_k(z),\quad i\in S_n^c,
$$
or
$$
\boldsymbol\psi(z)=z\bold p_n+zP_n\boldsymbol\psi(z),\quad \bold p_n:=
\bigl[\{p(i,A)\}_{i\in S_n^c}\bigr]^T.
$$
Therefore, introducing the $S_n^c\times S_n^c$ identity matrix $I_n$ and
using (3.16)-(3.17),
$$
\align
\boldsymbol\psi(z)=&\,z(I_n-zP_n)^{-1}\bold p_n=z\sum_{\nu\ge 0}z^{\nu}(P_n)^{\nu}
\bold p_n\\
=&\,\frac{z}{1-z\lambda_n}\,\bold f_n(\boldsymbol\sigma_n\bold p_n)+\bold R_n(z); \tag 3.18\\
\bold R_n(z):=&\,z\sum_{\nu\ge 0}z^{\nu}R_{n,\nu}\bold p_n
\endalign
$$
By (3.17) and
$$
\|\bold p_n\|_{L_{\infty}(S_n^c)}\le p(A):=\sup_{i\in \Bbb N}p(i,A),
$$
we have: each component of $\bold R_n(z)$ is analytic for $|z|<(1-\delta_0/3)^{-1}$, and
$$
\|\bold R_n(z)\|_{L_{\infty}(S_n^c)}\le \frac{2p(A)}{1-|z|(1-\delta_0/3)}.
$$
\si

Therefore each $\psi_i(z)$ initially defined in the unit disk admits a {\it meromorphic\/}
extension to the open disk of radius $(1-\delta_0/3)^{-1}>1$, with a single, simple pole
$z=1/\lambda_n$ in that disk.
\si

As for the explicit term in (3.18),
$$
\bigl[\bold f_n(\boldsymbol{\sigma}_n\bold p_n)\bigr]_i=C(A)f_n(i),
\quad C(A):=\sum_{j\in S_n^c}
\sigma_n(j)p(j,A).
$$

In particular, setting $z=1$,
$$
P_i\{X(T(S_n))\in A\}=\psi_i(1)=\frac{C(A)f_n(i)}{1-\lambda_n}+
O\bigl(p(A)\bigr).
$$
Since $P_i\{X(T(S_n))\in S_n\}=1$, we then obtain
$$
1-\lambda_n=\frac{C(S_n)f_n(i)}{1+\e_i(n)},\quad \e_i(n)=O(p(S_n)).\tag 3.19
$$
Therefore
$$
P_i\{X(T(S_n))\in A\}=\frac{C(A)}{C(S_n)} +O(p(S_n)),\quad i\in S_n^c.
$$
This uniform estimate and (2.23), with $U_n=A_n$, easily imply that
$$
\frac{C(A)}{C(S_n)}=\frac{\pi(A)}{\pi(S_n)}+O(p(S_n)). \tag 3.20
$$
Furthermore, given a positive integer $\tau$,
$$
P_i\{T(S_n)=\tau,\,X(T(S_n))\in A\}=[z^{\tau}]\psi_i(z)=\frac{1}{2\pi i}
\oint_{L}\frac{\psi_i(z)}{z^{\tau+1}}\,dz,
$$
where $L$ is a circular contour $|z|=1$. By (3.18), the {\it extended\/} $\psi_i(z)$ has a
unique singularity, a simple pole, in a ring between $L$ and $L_1$,  which is
the circular contour of radius $(1-\delta_0/4)^{-1}$. Using (3.18) and the residue theorem,
we obtain
$$
\align
\frac{1}{2\pi i}\oint_{L}\frac{\psi_i(z)}{z^{\tau+1}}\,dz=&\,
-\frac{C(A_n)f_n(i)}{2\pi i}\oint_{L_1}\frac{1}{(1-z\lambda_n)z^{\tau}}\,dz
+O\bigl((1-\delta_0/4)^{\tau}p(A)\bigr)\\
=&\,C(A_n) f_n(i)\lambda_n^{\tau-1}+O\bigl((1-\delta_0/4)^{\tau}p(A)\bigr)
\endalign
$$
Thus, by (3.19) and (3.20),
$$
\multline
P_i\{T(S_n)=\tau,\,X(T(S_n))\in A\}=
(1-\lambda_n)\lambda_n^{\tau-1}\frac{\pi(A)}{\pi(S_n)}\\
+O\bigl[(1-\lambda_n)\lambda_n^{\tau-1}p(S_n)\bigr]+O[(1-\delta_0/4)^{\tau}p(S_n))].
\endmultline
\tag 3.21
$$
In particular, for $A=S_n$,
$$
\multline
P_i\{T(S_n)=\tau\}=(1-\lambda_n)\lambda_n^{\tau-1}\\
+O\bigl[(1-\lambda_n)\lambda_n^{\tau-1}p(S_n)\bigr]+O[(1-\delta_0/4)^{\tau}p(S_n))].
\endmultline\tag 3.22
$$
Now we had proved already that, under the tightness only,
$$
E_i[T^k(S_n)]\sim k! E_i^k[T(S_n)]\sim k! \pi^{-k}(S_n),
$$
so that
$$
E_i[(T(S_n))_k]\sim k!\pi^{-k}(S_n).
$$
According to (3.22), we also have
$$
E_i[(T(S_n))_k]=\frac{k!}{(1-\lambda_n)^k}+O(p(S_n)(1-\lambda_n)^{-k}).
$$
Comparing the two formulas we see that $1-\lambda_n\sim\pi(S_n)$. In fact,
we can say more. From (3.22) it follows that, uniformly for $i\in S_n^c$,
$$
E_i[T(S_n)]=\frac{1+O(p(S_n))}{1-\lambda_n}.
$$
Combining this with (2.10), we get
$$
1-\lambda_n=\pi(S_n)\bigl(1+O(p(S_n))\bigr).\tag 3.23
$$
The rest is short. Let $C\subseteq S_n\times\Bbb N$, and 
$C_{\tau}=\{k\in S_n\,:\,(k,\tau)\in C\}$. From (3.21) it follows that,  uniformly for
$i\in S_n^c$ {\it and\/} $C$,
$$
P_i\{(X(T(S_n)),T(S_n))\in C\}=\sum_{\tau\in \Bbb N}(1-\lambda_n)\lambda_n^{\tau-1}
\frac{\pi(C_{\tau})}{\pi(S_n)}+O(p(S_n)). \tag 3.24
$$
And, by (3.23),
$$
\multline
\left|\sum_{\tau\in \Bbb N}(1-\lambda_n)\lambda_n^{\tau-1}\frac{\pi(C_{\tau})}{\pi(S_n)}-
\sum_{\tau\in \Bbb N}\pi(S_n)(1-\pi(S_n))^{\tau-1}\frac{\pi(C_{\tau})}{\pi(S_n)}\right|\\
\le\bigl|(1-\lambda_n)-\pi(S_n)\bigr|\sum_{\tau\in \Bbb N}x^{\tau-2}[1+\tau (1-x)]\,
\frac{\pi(C_{\tau})}{\pi(S_n)},
\quad (x\text{ between }
\lambda_n\text{ and }1-\pi(S_n))\\
\le\bigl|(1-\lambda_n)-\pi(S_n)\bigr|\left[2(1-x)\sum_{\tau\ge 1}\tau x^{\tau-1}+
2\sum_{\tau\ge 1}x^{\tau-1}\right]\\
=\bigl|(1-\lambda_n)-\pi(S_n)\bigr|\cdot 4(1-x)^{-1}
=O(p(S_n)).
\endmultline
$$
So (3.24) becomes 
$$
\|\partial_{n,i}-\partial_n\|_{TV}=O(p(S_n)),\quad i\in S_n^c.
$$

Suppose that $i\in S_n$. Then
$$
\multline
P_i\{(X(T(S_n)),T(S_n))\in C\}\\
=\,\sum_{k\in S_n^c}p(i,k)P_k\{(X(T(S_n)),T(S_n)+1)\in C\}+O(p(S_n)),
\endmultline\tag 3.25
$$
where, by (3.24),
$$
\multline
P_k\{(X(T(S_n)),T(S_n)+1)\in C\}=\sum_{\tau\in \Bbb N}(1-\lambda_n)\lambda_n^{\tau-1}
\frac{\pi(C_{\tau+1})}{\pi(S_n)}+O(p(S_n))\\
=\sum_{\tau\ge 2}(1-\lambda_n)\lambda_n^{\tau-2}\,
\frac{\pi(C_{\tau})}{\pi(S_n)}+O(p(S_n))=\sum_{\tau\ge 2}(1-\lambda_n)\lambda_n^{\tau-1}
\frac{\pi(C_{\tau})}{\pi(S_n)}+O(1-\lambda_n+p(S_n))\\
=\sum_{\tau\in\Bbb N}(1-\lambda_n)\lambda_n^{\tau-1}\,
\frac{\pi(C_{\tau})}{\pi(S_n)}+O(1-\lambda_n+ p(S_n))\\
=\sum_{\tau\in \Bbb N}(1-\lambda_n)\lambda_n^{\tau-1}
\frac{\pi(C_{\tau})}{\pi(S_n)}+O(p(S_n)).
\endmultline
$$
Therefore, by (3.25), (3.24) holds for $i\in S_n$ as well. This completes the proof of
Theorem 3.1.\qed
\bi

Let $T_{n,r}$ be the time intervals between
consecutive visits to $S_n$. So $T_{n,1}=T(S_n)$, and, for $r>1$, 
$$
\align
T_{n,r}=&\,\min\left\{t >\Cal T_{n,r-1}\,:\,X(t)\in S_n\right\}-\Cal T_{n,r-1},\\
\quad \Cal T_{n,r-1}:=&\,\sum_{k<r}T_{n,k},
\endalign
$$
i. e. $\Cal T_{n,r}$ is the time of $r$-th visit to $S_n$. Let $X_{n,r}=X(\Cal T_{n,r})$, i. e.
$X_{n,r}$ is a state in $S_n$ visited at time $\Cal T_{n,r}$. Introduce a random sequence
$\{\ell_r;t_r\}_{r\ge 1}$, where all $\ell_1,t_1,\ell_2,t_2,\dots$ are independent and, for each $r$,
$$
P\{\ell_r\in A\}=\frac{\pi(A)}{\pi(S_n)},\quad A\subseteq S_n,
$$
while $t_r$ is distributed geometrically, with success probability $\pi(S_n)$. Also, for two
random vectors, $\bold Y$ and $\bold Z$, of a common dimension $\nu$, let $d_{TV}(\bold Y,\bold
Z)$ denote the total variation distance between the distributions of $\bold Y$ and $\bold Z$, i. e.
$$
d_{TV}(\bold Y;\bold Z) =\sup_{B\in \Cal B^{\nu}} |P\{\bold Y\in B\}-P(\bold Z\in B)|.
$$
Since $|x|$ is convex,
$$
\multline
0.5\sup_{f:\|f\|_{ L_{\infty}(\Bbb N^{\nu})}\le 1}\bigl|E[f(\bold Y)]-E[f(\bold Z)]
\bigr|\le d_{TV}(\bold Y;\bold Z)\\
\le \sup_{f:\|f\|_{ L_{\infty}(\Bbb N^{\nu})}\le 1}\bigl|E[f(\bold Y)]-E[f(\bold Z)]
\bigr|.
\endmultline
$$

Theorem 3.1 implies the following.
\proclaim{Theorem 3.3} Uniformly for an initial state $i\in \Bbb N$, 
$$
d_{TV}\bigl(\{X_{n,r}\,;\, T_{n,r}\}_{1\le r\le k}\,;\,\{\ell_r;t_r\}_{1\le r\le k}\bigr)=O(kp(S_n)),
\tag 3.26
$$
Thus, if $k=k(n)$ is such that $kp(S_n)\to 0$, the random sequence
$\{X_{n,r}\, ;\mathbreak \,T_{n,r}\}_{1\le r\le k}$ is asymptotic,
with respect to the total variation distance, to  the \linebreak Bernoulli sequence 
$\{\ell_r,t_r\}_{1\le r\le k}$. 
\endproclaim
\bi
{\bf Proof of Theorem 3.3.\/} We prove (3.26) by induction on $k$. For $k=1$, it is 
the statement of Theorem 3.1.  Assume (3.26) holds for some $k\ge 1$. 
Let  $f:\Bbb N^{k+1}\times \Bbb N^{k+1}\to\Bbb R$ have  $
\|f\|_{L_{\infty}(\Bbb N^{k+1}\times \Bbb N^{k+1})}\le 1$. Denote
$$
\aligned
&\bold X=\{X_{n,r}\}_{1\le r\le k+1},\quad\bold X^{(k)}=\{X_{n,r}\}_{1\le r\le k}, \\
&\bold Y=\{T_{n,r}-\delta(i,r)\}_{1\le r\le k+1},\quad \bold Y^{(k)}=\{T_{n,r}-\delta(i,r)\}_{1\le r\le k},\\
&\bold x^{(k)}=\{x_r\}_{1\le r\le k},\quad\bold y^{(k)}=\{y_r\}_{1\le r\le k}.
\endaligned\tag 3.27
$$
We write first
$$
\multline
E_i[f(\bold X\,;\,\bold Y)]\\
=E_i\bigl[E_i[f\bigl((\bold X^{(k)},X_{n,k+1})\,;\,(\bold Y^{(k)},T_{n,k+1}-1)\bigr)\,|\,
(\bold X^{(k)},\bold Y^{(k)})]\bigr].
\endmultline \tag 3.28
$$
By the strong Markov property,  
$$
\multline
E_i[f\bigl((\bold X^{(k)},X_{n,k+1})\,;\,(\bold Y^{(k)},T_{n,k+1}-1)\bigr)\,|\,
(\bold X^{(k)},\bold Y^{(k)})]\\
=\left.E_{x_k}[f((\bold x^{(k)},X),
(\bold y^{(k)},T-1))\bigr]\right|_{\bold x^{(k)}=\bold X^{(k)},\,\bold y^{(k)}=\bold Y^{(k)}};
\endmultline\tag 3.29
$$
here  $X,T$ are the location and the time of the first hit of
$S_n$ for the chain starting at $x_k\in S_n$. Using (3.26) for $k=1$, we have
$$
\left|E_{x_k}[f((\bold x^{(k)},X),
(\bold y^{(k)},T-1))\bigr]-E[f((\bold x^{(k)},\ell),
(\bold y^{(k)},t))\bigr]\right| =O(p(S_n)),\tag 3.30
$$
uniformly for $j\in S_n$. (Here $(\ell,t)\overset{\Cal D}\to\equiv (\ell_r,t_r)$.) So, introducing 
$\tilde f:\Bbb N^k\times\Bbb N^k\to\Bbb R$ by
$$
\tilde f(\bold x^{(k)},\bold y^{(k)})=E\left[f\bigl((\bold x^{(k)},\ell),
(\bold y^{(k)},t)\bigr)\right],
$$
and using (3.28)-(3.30), we have
$$
\bigl|E_i[f(\bold X,\bold Y)]-E_i[\tilde f(\bold X^{(k)},\bold Y^{(k)})]\bigr| =O(p(S_n)).
 \tag 3.31
$$
Besides, applying the inductive hypothesis to $\tilde f$, we also have
$$
E_i[\tilde f(\bold X^{(k)},\bold Y^{(k)})]-E\bigl[\tilde f(\{\ell_r,t_r\}_{1\le r\le k})\bigr]
=O(kp(S_n)).\tag 3.32
$$
It follows from Fubini theorem and (3.31)-(3.32),  that
$$
\align
E_i\bigl[f(\bold X,\bold Y)\bigr]-E\bigl[f(\{\ell_r,t_r\}_{1\le r\le k+1}))\bigr]
=&\,E_i\bigl[f(\bold X,\bold Y)\bigr]-E\bigl[\tilde f(\{\ell_r,t_r\}_{1\le r\le k}))\bigr]\\
=&\,O((k+1)p(S_n)),
\endalign
$$
which proves the inductive step. So (3.26) holds for all $k$. \qed
\bi

Let us apply Theorem 3.3 to the extreme values for the t.e.m.~ chains. Given a large $N$, 
let $X^{(j)}=X^{(N,j)}$ denote the $j$-th largest among
$X(1),\dots,X(N)$; in particular, $X^{(1)}=\max_{1\le t\le N}X(t)$. From now on we 
will use a notation $S_n=\{n+1,n+2,\dots\}$.
\proclaim{Corollary 3.4} Uniformly for $i=X(0)$,
$$
P_i\{X^{(\mu)}\le n\} =\, P\bigl\{\text{Poisson }(N\pi(S_n))<\mu\bigr\}
+O\bigl(\mu^2/N+Np^2(S_n)\bigr).
\tag 3.33
$$
\endproclaim
\noindent
{\bf Proof of Corollary 3.4.\/} $X^{(\mu)}\le n$ iff during $[1,N]$ the chain visited $S_n$
at most $\mu-1$ times. So, by Theorem 3.3,
$$
P_i\{X^{(\mu)}\le n\}=\sum_{j<\mu}\binom{N}{j}\pi^j(S_n)(1-\pi(S_n))^{N-j}+
O(\mu p(S_n)).\tag 3.34
$$
Here
$$
(N)_j(1-\pi(S_n))^{-j}=N^j\bigl(1+O(\mu^2/N) +O(\mu p(S_n))\bigr),
$$
and
$$
\align
(1-\pi(S_n))^N- e^{-N\pi(S_n)}\le&\, 2Ne^{-N\pi(S_n)}
\bigl(e^{-\pi(S_n)}-(1-\pi(S_n))\bigr)\\
\le&N\pi^2(S_n)e^{-N\pi(S_n))}.
\endalign
$$
So, as $\pi(S_n)\le p(S_n)$, (3.34) becomes (3.33).\qed
\bi

Corollary (3.4) is a special case of the following result. Given $a<b\le\infty$,
denote $S_{a,b}=(a,b]$, i. e. $S_{a,b}=S_a\setminus S_b$. Let $V_{a,b}=V_{N,a,b}$
denote the number of visits to $S_{a,b}$  during $[1,N]$,  and $\lambda_{a,b}=
\lambda_{N,a,b}=N\pi(S(a,b))$.
\proclaim{Theorem 3.5} Let $(a_1,b_1],\dots (a_k,b_k]$ be disjoint. Uniformly for $i=X(0)$,
$$
\aligned
P_i\left\{\bigcap_{1\le\ell\le k}\bigl\{V_{a_\ell,b_{\ell}}\le \mu_{\ell}\bigr\}
\right\}=&\,\prod_{1\le\ell\le k}P\bigl\{\text{Poisson }(\lambda_{a_{\ell},b_{\ell}})\le\mu_{\ell}
\bigr\}\\
&+O\bigl(\mu^2/N+Np^2(S_{a})\bigr),
\endaligned\tag 3.35
$$
where $\mu=\mu_1+\cdots+\mu_k$, $a=\min_{\ell}\,a_{\ell}$. Thus, if $\mu=o(N^{1/2})$ and 
$Np^2(S_a)=o(1)$, the numbers of visits to non-overlapping intervals $(a_{\ell},b_{\ell}]$ are
asymptotically independent Poissons with parameters $N\pi(S_{a_{\ell},b_{\ell}})$.
\endproclaim
\noindent
{\bf Proof of Theorem 3.5.\/} Applying Theorem 3.3 to 
$S:=\bigcup\limits_{1\le\ell\le k}S_{a_{\ell},b_{\ell}}$,
$$
\align
P_i\left\{\bigcap_{1\le\ell\le k}\bigl\{V_{a_\ell,b_{\ell}}\le \mu_{\ell}\bigr\}
\right\}
=&\,\sum_{j_1\le\mu_1;\dots;j_k\le\mu_k}\binom{N}{j_1,\dots,j_k}
\prod_{1\le\ell\le k}\pi^{j_{\ell}}(S_{a_{\ell},b_{\ell}})\\
&\times \left(1-\sum_{1\le\ell\le k}\pi(S_{a_{\ell},b_{\ell}})\right)^
{N-j}+O(\mu p(S_a)),
\endalign
$$
where $j=j_1+\cdots+j_{\ell}$, and
$$
\binom{N}{j_1,\dots,j_k}=\frac{N!}{j_1!\cdots j_k!\, (N-j)!}.
$$ 
The rest runs parallel with the proof of
Corollarry 3.4. \qed
\bi

Analogously we obtain a relatively simple asymptotic formula for the joint distribution of
$X^{(1)},\dots,X^{(\mu)}$.
\proclaim{Theorem 3.6} Let $\infty=n_0\ge n_1\ge n_2\ge\cdots\ge n_{\mu}$. Uniformly for
$i=X(0)$,
$$
\multline
P_i\left\{\bigcap_{1\le\ell\le\mu}\bigl\{X^{(\ell)}\le n_{\ell}\bigr\}\right\}\\
=
\sum_{\nu_1,\dots,\nu_{\mu}\atop
\forall r\le\mu\,:\, \sum_{j=1}^r\nu_j\le r-1}\prod_{1\le r\le \mu}\!\!\!\!
P\bigl\{\text{Poisson }(\lambda_{n_r, n_{r-1}})= \nu_r\bigr\} +O\bigl(\mu^2/N+
Np^2(S_{n_{\mu}})
\bigr).
\endmultline\tag 3.36
$$
More generally, let
$$
B\subseteq\bigl\{\bold x=(x_1,\dots,x_{\mu})\in \Bbb N\,:\,
x_1\ge\cdots\ge x_{\mu}\bigr\}.
$$
Given $\bold x$, let $y_1(\bold x)>\cdots>y_m(\bold x)$ denote all the
distinct values (range) of the sequence $x_1,\dots,x_{\mu}$, and let $a_j=a_j
(\bold x)>0$
be the multiplicity of $y_j=y_j(\bold x)$. So $m=m(\bold x)\le \mu$, and 
$a_1+\cdots +a_m=\mu$.
Then, denoting $n(B)=\inf\limits_{\bold x\in B}x_{\mu}$,  and setting $y_0=\infty$,
$$
\multline
P_i\bigl\{(X^{(1)},\dots,X^{(\mu)})\in B\bigr\}\\
=\sum_{\bold x\in B}\,\,\prod_{1\le r\le m}e^{-N\pi([y_r,y_{r-1}))}
\,\frac{(N\pi(y_r))^{a_r}}{a_r!}
+O\bigl(\mu^2/N +Np^2(S_{n(B)})\bigr)\\
=\sum_{\bold x\in B}\,\,e^{-N\pi([y_{m},\infty))}(N\pi([y_{m},\infty))^{\mu}
\prod_{1\le r\le m}
\,\frac{\sigma^{a_r}(y_r)}{a_r!}
+O\bigl(\mu^2/N +Np^2(S_{n(B)})\bigr),
\endmultline\tag 3.37
$$
where $\sigma(y)=\pi(y)/\pi([y_{m},\infty))$,  $y\in [y_{m},\infty)$.
\endproclaim
\bi

In the next section we will describe two models of a random constrained composition, and show
that each random composition is sharply approximated by a t.e.m.~ chain. It will
enable us to use Corollary 3.4 and Theorems 3.5, 3.6 for analysis of the limiting
distribution of the larger parts.
\bi
{\bf 4. Two random constrained compositions and Markov chain approximations.\/} We focus on
two interesting cases of such compositions, the column-convex-animals (cca)
compositions and the Carlitz (C) compositions.
\bi
{\bf 4.1. Defintions and some basic facts.\/} {\bf (a)\/} A column-convex animal (cca) is a sequence 
of contiguous vertical segments of unit squares in $\Bbb Z^2$, 
ordered from left
to right, such that every two successive columns
have a common boundary consisting of at least one vertical edge of $\Bbb Z^2$. If the total
number of unit squares involved is $\nu$ then the lengths of the vertical segments form
a composition of $\nu$; we call it a cca composition. Let $T(\nu,\mu)$ denote 
the 
total number
of the cca compositions of $\nu$ with $\mu$ parts; then $T(\nu):=\sum_{\mu\ge 1}
T(\nu,\mu)$ is
the total number of the cca compositions of $\nu$. Introduce $f(w,z)$, 
the bivariate generating 
function
(BGF) of $T(\nu,\mu)$,
$$
f(w,z)=\sum_{\mu,\,\nu\ge 1}T(\nu,\mu) w^{\mu}z^{\nu}.
$$
Louchard [19] found that
$$
\aligned
f(w,z)=&\frac{wz(z-1)^3}{h(w,z)},\\
h(w,z):=&\,z^4(w-1)+z^3(w^2-w+4)-z^2(w+6)+z(w+4)-1.
\endaligned\tag 4.1.1
$$
\bi
Therefore $f(z)$, the GF of $T(\nu)$, is
$$
f(z)=f(1,z)=\frac{z(z-1)^3}{h(1,z)}=\frac{(z-1)^3}{4z^3-7z^2+5z-1},\tag 4.1.2
$$
a formula discovered earlier by Klarner [15]. Privman and Forgacs [22] used (4.1.2),
and Darboux theorem, to show
that
$$
T(\nu)=\frac{C}{z_*^{\nu}}\bigl(1+O(\gamma^{\nu})\bigr),\tag 4.1.3
$$
where $\gamma<1$,  $C=0.18\dots$, and $z_*=0.31\dots$ is the smallest-modulus solution
of $h(1,z)=0$. 
\si

We get a uniformly random  cca composition of $\nu$, if we assume that each  composition 
has the same probability, $1/T(\nu)$.  It was discovered in [19], [20] that
the distribution of the last (first) part is asymptotic to
$$
\aligned
\pi_1(k)=&\,z_*^k(k+a),\\
a:=&\,\frac{1-z_*}{z_*}-\frac{1}{1-z_*}=0.75\cdots,
\endaligned\tag 4.1.4
$$
which is directly seen as a probability distribution.  Besides, the joint
distribution of {\it two\/} consecutive parts $Y_t$ and $Y_{t+1}$, with both $t$ and 
$\Cal M-t$ of order $\Theta(\nu)$, was shown to be asymptotic to that of two consecutive states of an ergodic
Markov chain on $\Bbb N$, in a stationary regime, with transition probabilities
$$
p(i,k)=z_*^k(i+k-1)\frac{k+a}{i+a},\tag 4.1.5
$$
and a stationary distribution
$$
\aligned
\pi(k)= &\,A^{-1} z_*^k(k+a)^2,\\
 A:=&\,\sum_{k\ge 1}z_*^k(k+a)^2=
\frac{z_*^2}{(1-z_*)^3}+\frac{1-z_*}{z_*}. 
\endaligned\tag 4.1.6
$$
That $\sum_{k\ge 1}p(i,k)=1$ follows from another formula for $a$,
$$
a=\frac{2z_*^2}{(1-2z_*)(1-z_*)}.
$$
(The given formulation is slightly different from, but equivalent to that in [19], 
[20].)
One way to derive (4.14) is to use (4.1.3) and a formula for $f_k(z)$, the generating function of
the cca compositions with the first (last) part equal $k$,
$$
f_k(z)=z^k+z^k f(z)\left[k+\frac{z^3-z^2+z}{(1-z)^3}\right],\tag 4.1.7
$$
which can be read out of [19]. Comparing the first line in (4.1.4) and (4.1.7) we must also have
yet another formula for $a$, namely
$$
a=\frac{z_*^3-z_*^2+z_*}{(1-z_*)^3},\tag 4.1.8
$$
which is indeed the case.
\bi

{\bf (b)\/} A Carlitz (C) composition of $\nu$ is defined as a composition  such that every 
two consecutive parts are distinct from each other. The counterparts of the cited
results for the cca compositions are as follows.  Carlitz [8] proved that
$$
\aligned
f(w,z)=&\,-1+\frac{1}{h(w,z)},\\
h(w,z):=&\,1-\sum_{j\ge 1}(-1)^{j+1}\frac{w^jz^j}{1-z^j};
\endaligned\tag 4.1.9
$$
for $|w|\le 1$, $h$, as a function of $z$ is analytic for $|z|<1$, and for $|w|\ge 1$, $h$ is analytic
for $|z|<1/|w|$. Louchard and Prodinger [21] found a rather more tractable expression
for $h$, namely
$$
h(w,z)=1-\sum_{j\ge 1}\frac{wz^j}{1+wz^j}.\tag 4.1.10
$$
(4.1.9) and (4.1.10)  were used in [20]  to show that
$$
T(\nu)=\frac{C}{z_*^{\nu}}\big(1+O(\gamma^{\nu})\bigr),\tag 4.1.11
$$
where $\gamma<1$, $C=0.456\dots$, and $z_*=0.57\dots$ is the smallest-modulus solution
of $h(1,z)=0$.
\si

We get a uniformly random C-composition of $\nu$, if we assume that each C-composition 
has the same probability, $1/T(\nu)$. In a striking analogy with the random cca composition, the two consecutive
parts $Y_{t}$ and $Y_{t+1}$, deep inside the composition, are also jointly asymptotic to 
the two consecutive states
of an ergodic Markov chain, with transition probabilities
$$
p(i,k)=\left\{\alignedat2
&z_*^k\,\frac{1+z_*^ i}{1+z_*^k},\quad&&i\neq k,\\
&0,\quad&& i=k.\endalignedat\right.\tag 4.1.12
$$
and a stationary distribution
$$
\pi(k)=A^{-1}\frac{z_*^k}{(1+z_*^k)^2},\quad A:=\sum_{k\ge 1}\frac{z_*^k}{(1+z_*^k)^2}.
\tag 4.1.13
$$
And the limiting distribution of $Y_1$ is
$$
\pi_1(k)=\frac{z_*^k}{1+z_*^k},\tag 4.1.14
$$
which follows from (4.1.11) and a counterpart of (4.1.7),
$$
f_k(z)=\frac{z^{k+1}}{1+z^{k+1}}+f(z)\,\frac{z^k}{1+z^k}.\tag 4.1.15
$$
(That (4.1.12) and (4.1.14) and  are indeed probability distributions follows 
from the definition 
of $z_*$ as a root of $h(1,z)=0$ and (4.1.10).)
\bi

For each of the compositions, an equation $h(w,z)=0$ (for the attendant function $h(w,z)$)
determines a root $z(w)$,
well defined for $w$ sufficiently close to $1$, such that $z(1)=z_*$,
$z(w)$ is infinitely differentiable, and $z^\prime(1)<0$. The number of parts $\Cal M$
for each of the random compositions was shown, in [19] and [21] resp.,
to be Gaussian in the limit $\nu\to\infty$, with mean $\alpha\nu$ and
variance $\beta\nu$, where
$$
\alpha=-\frac{z^\prime(1)}{z(1)}=-\frac{z^\prime(1)}{z_*},
\quad \beta=\alpha^2+\alpha-\frac{z^{\prime\prime}(1)}{z_*}.\tag 4.1.16
$$
In particular,
$$
\alpha=\left\{\alignedat2
&-\frac{12z_*^2-14z_*+5}{z_*^4+z_*^3-z_*^2+z_*}=0.45\dots,\quad&& (\text{for cca}),\\
&-\frac{\sum_{j\ge 1}\frac{jz_*^{j-1}}{(1+z_*^j)^2}}
{\sum_{j\ge 1}\frac{z_*^j}{(1+z_*^j)^2}}=0.35\dots,\quad&& (\text{for C});
\endalignedat\right.\tag 4.1.17
$$
needless to say, in each case $z_*$ is the root of the corresponding
equation $h(1,z)=0$.
\si

In Appendix we will prove the following large deviation result.
\proclaim{Lemma 4.1.1} For each of the compositions, there exists an absolute constant 
$c>0$ such that 
$$
P\{|\Cal M-\alpha\nu |\ge s\}\le c\nu \exp(-s^2/3\beta\nu),
$$
provided that $s=o(\nu)$. Thus 
$$
P\{|\Cal M-\alpha\nu|\le \nu^{1/2}\ln\nu\}\ge 1 - \nu^{-K},\quad\forall\, K>0.\tag
4.1.18
$$ 
\endproclaim
\si
{\bf Note.\/} Borrowing a term from Knuth et al. [17], the event on
the left of (4.1.18) happens {\it quite surely\/} (q.s.).
\bi 
{\bf 4.2. Approximating the random compositions by the Markov chains.\/} The 
results cited above strongly suggest, though not actually prove, that the random cca 
composition and
the random C-composition considered as random processes are each asymptotic to 
its own Markov chain, defined in (4.1.4)-(4.1.5) and (4.1.11)-(4.1.13) respectively.
\si

The following theorem confirms this natural conjecture.
\proclaim{Theorem 4.2.1} Let $\bold Y=\{Y_t\}_{t\ge 1}$ be either the random cca 
composition, or the random 
C-composition of $\nu$. Let $\bold Z=\{Z(t)\}_{t\ge 1}$ be the corresponding Markov chain
with the transition probabilities $p(i,k)$, and $Z(1)$ having the distribution $\{\pi_1(i)\}_{i\ge 1}$.
Introduce 
$$
\aligned
\hat\Cal M=&\,\max\left\{1\le m<\Cal M\,:\,Y_1+\cdots+Y_m\le \nu-\ln^2 \nu\right\},\\
 \hat M=&\,
\max\left\{m\ge 1\,:\,Z(1)+\cdots+Z(m)\le \nu-\ln^2 \nu\right\};
\endaligned\tag 4.2.1
$$
in particular, $\hat\Cal M\in (\Cal M-\ln^2\nu-1,\Cal M)$. Let $\hat\partial$ and $\hat d$
denote the probability distribution of $(\hat\Cal M,\,(Y_1,\dots,Y_{\hat\Cal M}))$
and $(\hat M,\,(Z(1),\dots,Z(\hat M))$ respectively. 
For each chain, 
$$
\|\hat\partial-\hat d\|_{TV}=O\bigl(\nu^{-K}\bigr),\quad\forall\,K>0.\tag 4.2.2
$$
\endproclaim
\si

So, the random composition of $\nu$, read from left to right, is closely
approximated by the corresponding Markov chain, as long as the accumulated
sum of parts stays below $\nu-\ln^2\nu$. (A restriction of this sort is
unavoidable: like the first part, the last part of the random composition has 
the distribution $\pi_1$,
which differs from the stationary distribution $\pi$.) Now, we will see that, with high probability,
the extreme-valued parts are in this ``bulk'' of the composition, implying that
they are well approximated by the extreme-valued states of the $\hat M$-long
segment of the corresponding Markov chain. It is easy to verify that
$$
\sup_i\sum_{k>n}p(i,k)=\left\{\alignedat2
&O(z_*^n n^2),\quad&&\text{cca chain},\\
&O(z_*^n),\quad&&\text{C-chain},\endalignedat\right.\tag 4.2.3
$$
where $z_*=0.31\dots$ for the cca chain and $z_*=0.57\dots$ for the C-chain. That is,
the chains meet the tightness condition (2.1).  And the exponential mixing 
property in the form of (3.3) is easily verified as well. So we are able to use 
Corollary 3.4 and Theorem 3.6, say, for
derivation of the limiting distribution of those extreme values, and then the
last theorem for a quick proof of the corresponding results regarding extreme
-valued parts of each of the random compositions.

Turning the tables, we can also use 
Theorem 4.2.1 and Lemma 4.1.1 to determine the
very likely bounds of $\hat M$ 
with sufficient
accuracy. Since $\hat\Cal M\in (\Cal M-\ln^2\nu-1,\Cal M)$, Lemma 4.1.1 implies 
that q.s.
$$
|\hat\Cal M-\alpha\nu|\le 2\nu^{1/2}\ln\nu.
$$
So, applying Theorem 4.2.1, we immediately see that
$$
|\hat M-\alpha\nu|\le 2\nu^{1/2}\ln\nu\tag 4.2.4
$$
q.s. as well. (!) 
\bi  
{\bf Proof of Theorem 4.2.1.\/} The key element is the following claim.
\si
\proclaim{Lemma 4.2.2} Let $\bold Y$ be either the random cca composition or the random 
C-composition of $\nu$. Let $k\ge 1$, $\bold i=(i_1,\dots,i_k)\in \Bbb N^k$, where $i_1+
\cdots+i_k<\nu$. Denote $P_{\nu}(\bold i)=P\{Y_1=i_1,\dots,Y_k=i_k\}$ and 
$P(\bold i)=P\{Z(1)=i_1,\dots, Z(k)=i_k\}$. Then, uniformly for $k$ and $\bold i$,
$$
P_{\nu}(\bold i)=P(\bold i)\exp\left(O\bigl(k\,\gamma\,^{\nu-|\bold i|}\bigr)\right),\quad
|\bold i|=i_1+\cdots+i_k,\tag 4.2.5
$$
where $\gamma$ comes from either (4.1.3) or (4.1.10).
\endproclaim
\si
{\bf Proof of Lemma 4.2.2.\/} Let $\bold Y$ be the random cca composition of $\nu$.
We will prove (4.2.5) by induction on $k$. 
\si

For $k=1$,
$$
P_{\nu}(i_1)=\frac{[z^{\nu}]\,f_{i_1}(z)}{[z^{\nu}]\,f(z)},\tag 4.2.6
$$
where $f(z)$ and $f_{i_1}(z)$ are given by (4.1.2) and (4.1.7) respectively.
Here, by (4.1.3),
$$
[z^{\nu}]\,f(z)=T(\nu)=\frac{C}{z_*^\nu}\,
\exp\left(O\bigl(\gamma\,^{\nu}\bigr)\right).\tag 4.2.7
$$
Further, by (4.1.7),
$$
\align
&[z^{\nu}]\,f_{i_1}(z)=\,\delta_{\nu,i_1} + [z^{\nu-i_1}]\,f(z)
\left[i_1+\frac{z^3-z^2+z}{(1-z)^3}\right]\\
=&\,i_1 T(\nu-i_1)+[z^{\nu-i_1}]\,f(z)\,
\frac{z^3-z^2+z}{(1-z)^3}\\
=&\,i_1\,\frac{C}{z_*^{\nu-i_1}}
\exp\left(O\bigl(\gamma\,^{\nu-i_1}\bigr)\right)
+\frac{C}{z_*^{\nu-i_1}}\,\frac{z_*^3-z_*^2+z_*}{(1-z^*)^3}
\exp\left(O\bigl(\gamma\,^{\nu-i_1}\bigr)\right).\tag 4.2.8
\endalign
$$
($z_*$ is the smallest modulus pole of $f(z) (z^3-z^2+z)(1-z)^{-3}$, as well.)
It follows from (4.2.6)-(4.2.8) and (4.1.8) that 
$$
\align
P_{\nu}(i_1)=&\,z_*^{i_1}\left[i_1+\frac{z_*^3-z_*^2+z_*}{(1-z_*)^3}\right]
\exp\left(O\bigl(\gamma\,^{\nu-i_1}\bigr)\right)\\
=&\,z_*^{i_1}(i_1+a)\exp\left(O\bigl(\gamma\,^{\nu-i_1}\bigr)\right)=
P(i_1)\exp\left(O\bigl(\gamma\,^{\nu-i_1}\bigr)\right),
\endalign
$$
which is (4.2.5) for $k=1$.
\si

Suppose that (4.2.5) holds for some $k\ge 1$. Let $\bold i=(i_1,\dots,i_{k+1})$
be such that $|\bold i|<\nu$. Let $\bold i^\prime=(i_2,\dots,i_{k+1})$; then
$|\bold i^\prime|<\nu-i_1$. Let $T(\bold i,\nu)$ and $T(\bold i^\prime,\nu-i_1)$
denote the total number of the cca of area $\nu$ ($\nu-i_1$ resp.) with
the first $k+1$ parts $i_1,\dots,i_{k+1}$ (the first $k$ parts $i_2,\dots,i_{k+1}$
resp.). By the definition of the cca composition,
$$
T(\bold i,\nu)=(i_1+i_2-1)T(\bold i^\prime,\nu-i_1).
$$
Therefore
$$
\align
P_{\nu}(\bold i)=&\,\frac{T(\bold i,\nu)}{T(\nu)}\\
=&\frac{[z^{\nu}]\,f_{i_1}(z)}{T(\nu)}\cdot \frac{T(\nu-i_1)}{[z^{\nu}]\,f_{i_1}(z)}
\cdot\frac{(i_1+i_2-1)T(\bold i^\prime,\nu-i_1)}{T(\nu-i_1)}\\
=&\,P_{\nu}(i_1)\,\frac{\exp\left(O\bigl(\gamma\,^{\nu-i_1}\bigr)\right)}
{i_1+a}\cdot (i_1+i_2-1)P_{\nu-i_1}(\bold i^\prime)\,\\
=&\,P(i_1)\,\frac{\exp\left(O\bigl(\gamma\,^{\nu-i_1}\bigr)\right)}
{i_1+a}\cdot (i_1+i_2-1)P(\bold i^\prime)
\exp\left(O\bigl(k\,\gamma\,^{\nu-i_1-|\bold i^\prime|}\bigr)\right)\\
=&\,P(i_1)\,\frac{i_1+i_2-1}{i_1+a}\,P(\bold i^\prime)
\exp\left(O\bigl((k+1)\,\gamma\,^{\nu-|\bold i|}\bigr)\right),
\endalign
$$
and we observe that
$$
\align
\frac{i_1+i_2-1}{i_1+a}\,P(\bold i^\prime)=&\,
\frac{i_1+i_2-1}{i_1+a}\,P(i_2)\prod_{r=2}^kp(i_r,i_{r+1})\\
=&p(i_1,i_2)\prod_{r=2}^kp(i_r,i_{r+1})=\prod_{r=1}^kp(i_r,i_{r+1}).
\endalign
$$
Hence
$$
\align
P_{\nu}(\bold i)=&\,P(i_1)\prod_{r=1}^kp(i_r,i_{r+1})
\exp\left(O\bigl((k+1)\gamma\,^{\nu-|\bold i|}\bigr)\right)\\
=&\,P(\bold i)\exp\left(O\bigl((k+1)\gamma\,^{\nu-|\bold i|}\bigr)\right),
\endalign
$$
which completes the inductive proof of (4.2.5) for the random cca composition.
The proof for the random C-composition is similar, and we omit it. \qed
\bi

Lemma 4.2.2 implies the bound (4.2.2) of Theorem 4.2.1 without much
difficulty. Consider, for instance, the random cca composition of $\nu$. Let $m$, 
$\bold i=(i_1,\dots,i_m)$ be given. Clearly
$$
\multline
P\{\hat\Cal M\ge m,\,Y_1=i_1,\dots,Y_m=i_m\}\\
=P\{\hat M\ge m,\,Z(1)=i_1,\dots,Z(m)=i_m\}=0,
\endmultline
$$
unless $|\bold i|\le \nu-\ln^2 \nu$. In the latter case $m\le \nu-\ln^2\nu$,
and, by Lemma 4.2.2,
$$
\align
&P\{\hat\Cal M\ge m,\,Y_1=i_1,\dots,Y_m=i_m\}=P(\bold i)
\exp\left(O\bigl(m\,\gamma\,^{\nu-|\bold i|}\bigr)\right)\\
=&P\{\hat M\ge m,\,Z(1)=i_1,\dots,Z(m)=i_m\}
\exp\left(O\bigl(\nu\,\gamma^{\ln^2\nu}\bigr)\right),
\endalign
$$
uniformly for $m$ and $\bold i$ in question. Consequently, uniformly for all
$m$ and $B\subseteq \Bbb N^m$,
$$
\align
&P\{\hat\Cal M\ge m,\, (Y_1,\dots,Y_m)\in B\}\\
=&P\{\hat M\ge m,\,(Z(1),\dots,Z(m)\in B\}
\exp\left(O\bigl(\gamma^{0.5\ln^2\nu}\bigr)\right),
\endalign
$$
whence
$$
\aligned
&P\{\hat\Cal M= m,\, (Y_1,\dots,Y_m)\in B\}\\
=&P\{\hat M= m,\,(Z(1),\dots,Z(m))\in B\}+O\bigl(\gamma^{0.5\ln^2\nu}\bigr).
\endaligned\tag 4.2.9
$$
Let $D\subseteq \Bbb N^{\nu+1}$ be given. For $\bold z\in \Bbb N^k$, $k\le\nu+1$,
 we write $\bold z\in D$
if $\bold z$ is a projection of a point in $D$ on the first $k$ coordinates.
Noticing that $\hat\Cal M\le\nu$ and $\hat M\le\nu$, we obtain from (4.2.8): uniformly
for all $D\in \Bbb N^{\nu+1}$,
$$
\align
&P\left\{\bigl(\hat\Cal M,\, (Y_1,\dots,Y_{\hat\Cal M})\bigr)\in D\right\}\\
=&P\left\{\bigl(\hat M,\,(Z(1),\dots,Z(\hat M))\bigr)\in D\right\}+
O\bigl(\gamma^{0.5\ln^2\nu}\bigr).
\endalign
$$
This completes the proof of Theorem 4.2.1.\qed
\bi
{\bf 4.3. Limiting distributions of the extreme parts of the random compositions.\/} By (4.2.3), for each
of the two chains,  q.s.
$$
N_1+1\le \hat M\le 1+N_2,\quad N_{1,2}=\bigl\lfloor\alpha\nu\pm 2\nu^{1/2}\ln\nu\bigr\rfloor.
\tag 4.3.1
$$
So q.s. the extreme values of $\{Z(t)\}_{0<t\le\hat M}$ are sandwiched between those of 
$\{Z(t)\}_{0<t\le N_1+1}$ and $\{Z(t)\}_{0<t\le N_2+1}$. Picking a generic $N\in [N_1,N_2]$,
introduce $\{X(t)\}_{0\le t\le N}=\{Z(t)\}_{0<t\le N+1}$. Here $X(0)$ has distribution 
$\pi_1(\cdot)$.
 \bi
 
Let $X^{(\mu)}$ be the $\mu$-th largest among $X(t)$, $t\in [1,N]$, for $X(0)=i$, 
$i\in \Bbb N$. By
Corollary 3.4,
$$
P_i\{X^{(\mu)}\le n\}=P\bigl\{\text{Poisson }(N\pi(S_n))<\mu\bigr\}+
O\bigl(\mu^2/N+Np^2(S_n)\bigr),\tag 4.3.2
$$
where
$$
\pi(S_n)=\sum_{k>n}\pi(k),\quad p(S_n)=\sup_i\sum_{k>n}p(i,k).
$$
Here $p(S_n)=O(n^2z_*^n)$ for the cca chain, and $p(S_n)=O(z_*^n)$ for
the C-chain, see (4.2.3). (Again, $z_*=0.31\dots$ for the
cca chain, and $z_*=0.57\dots$ for the C-chain.) Turn to $\pi(S_n)$. For the cca chain, 
by (4.1.6),
$$
\pi(S_n)=n^2z_*^nB\bigl(1+O(n^{-1})\bigr),\quad B:=\frac{z_*^2(1-z_*)^2}{z_*^3+(1-z_*)^4}.
\tag 4.3.3
$$
For the C-chain, by (4.1.13),
$$
\pi(S_n)=Bz_*^n\bigl(1+O(z_*^n)\bigr),\quad B:=A^{-1}\frac{z_*}{1-z_*}.\tag 4.3.4
$$
\proclaim{Lemma 4.3.1 (cca chain)} Suppose that
$$
n=\frac{\ln\bigl[\lambda^{-1}BN(\ln N/\ln z_*)^2\bigr]}{\ln (1/z_*)}\in \Bbb N,\tag 4.3.5
$$
where $\lambda=o(\ln N)$. If $\mu=o(\ln N)$, then, uniformly for 
$i\in \Bbb N$,
$$
P_i\{X^{(\mu)}\le n\}=P\bigl\{\text{Poisson }(\lambda)<\mu\bigr\}
+O\bigl[(\lambda+\mu)/\ln N\bigr].\tag 4.3.6
$$
Equivalently, define $W_{N,\mu}$ by
$$
X^{(\mu)}=\frac{\ln\bigl[W_{N,\mu}^{-1}\,BN(\ln N/\ln z_*)^2\bigr]}{\ln (1/z_*)};\tag 4.3.7
$$
then, for $s=o(\ln N)$ such that 
$$
\frac{\ln\bigl[s^{-1}BN(\ln N/\ln z_*)^2]}{\ln (1/z_*)}\in \Bbb N,\tag 4.3.8
$$
we have
$$
P_i\{W_{N,\mu}\ge s\}=P\{W_{\mu}\ge s\}+
O\bigl[(s+\mu)/\ln N\bigr];\tag 4.3.9
$$
here $W_{\mu}=V_1+\cdots+V_{\mu}$, and $V_1,\dots,V_{\mu}$ are independent 
exponentials with unit mean.
\endproclaim
This Lemma implies the following cruder result. (We use a symbol $O_p(1)$ to denote
a random variable bounded in probability as $N\to\infty$.)
\proclaim{Corollary 4.3.2} If $\mu=o(\ln N)$, then, uniformly for $i\in \Bbb N$,
$$
X^{(\mu)}=\frac{\ln\bigl(\mu^{-1}N\ln^2 N\bigr)}{\ln(1/z^*)}+O_p(1).\tag 4.3.10
$$
\endproclaim
Here are the counterparts for the chain associated with the random C-compo
\linebreak sition.
\proclaim{Lemma 4.3.3 (C-chain)} Suppose that
$$
n=\frac{\ln(\lambda^{-1}BN)}{\ln(1/z_*)}\in \Bbb N,\tag 4.3.11
$$
where $\lambda=o(N^{1/2})$. If $\mu=o(N^{1/2})$, then, uniformly for $i\in \Bbb N$,
$$
P_i\{X^{(\mu)}\le n\}=P\bigl\{\text{Poisson }(\lambda)<\mu\bigr\}
+O\bigl[(\lambda^2+\mu^2)/N\bigr].\tag 4.3.12
$$
Equivalently, define $W_{N,\mu}$ by
$$
X^{(\mu)}=\frac{\ln(W_{N,\mu}^{-1}BN)}{\ln(1/z_*)};\tag 4.3.13
$$
then, for $s=o(N^{1/2})$ such that 
$$
\frac{\ln(s^{-1}BN)}{\ln(1/z_*)}\in\Bbb N,\tag 4.3.14
$$
we have
$$
P_i\{W_{N,\mu}\ge s\}=P\{W_{\mu}\ge s\}+O\bigl[(s^2+\mu^2)/N\bigr].
$$
\endproclaim
\si
\proclaim{Corollary 4.3.4 (C-chain)} If $\mu=o(N^{1/2})$, then, uniformly for $i\in\Bbb N$,
$$
X^{(\mu)}=\frac{\ln(\mu^{-1} N)}{\ln(1/z^*)}+O_p(1).\tag 4.3.15
$$
\endproclaim
\bi
{\bf Proof of Lemma 4.3.1 and Corollary 4.3.2.\/} {\bf (a)\/}  By (4.3.3), (4.3.5) and 
(4.2.3),
$$
N\pi(S_n)=\lambda+O\bigl(\lambda/\ln N\bigr),\qquad Np^2(S_n)=
O\bigl[N^{-1}(N\pi(S_n))^2\bigr]=O(\lambda^2/N).
$$
Then, for $j\le\mu$,
$$
(N\pi(S_n))^j=\lambda^j\bigl(1+O(\mu/\ln N)\bigr).
$$
So, by Corollary 3.4, (3.33), and (4.2.3), (4.3.3),
$$
\align
P_i\{X^{(\mu)}\le n\}=&\,\sum_{j<\mu}e^{-N\pi(S_n)}\frac{(N\pi(S_n))^j}{j!}+
O\bigl((\mu^2+\lambda^2)/N\bigr)\\
=&\,\sum_{j<\mu}e^{-\lambda}\frac{\lambda^j}{j!}+O\bigl((\lambda+\mu)/\ln N\bigr).
\endalign
$$
\si
{\bf (b)\/} Given $s>0$,
$$
\left\lfloor\frac{\ln\bigl[s^{-1}BN(\ln n/\ln z_*)^2\bigr]}{\ln(1/z_*)}\right\rfloor=
\frac{\ln\bigl[s_1^{-1}BN(\ln n/\ln z_*)^2\bigr]}{\ln(1/z_*)},
$$
where $s_1\in [s,\,sz_*^{-1})$. Using the definition of $W_{N,\mu}$ in (4.3.7) and the
asymptotic formula (4.3.9) we obtain then: for $s=o(\ln N)$,
$$
\multline
P\{z_*W_{\mu}\ge s\}+O\bigl[(s+\mu)/\ln N\bigr]\le P_i\{W_{N,\mu}\ge s\}\\
\le P\{W_{\mu}\ge s\}+O\bigl[(s+\mu)/\ln N\bigr].
\endmultline\tag 4.3.16
$$
For $\mu$ fixed, (4.3.16) implies that 
$$
\lim_{A\to\infty}\liminf_{N\to\infty}P_i\{W_{N,\mu}\in [A^{-1},A]\}=1,
$$
i. e. , in probability, $W_{N,\mu}$ is bounded away from zero and infinity, whence
$\ln W_{N,\mu}=O_p(1)$. Suppose $\mu\to\infty$. Then $(W_{\mu}-\mu)/\mu^{1/2}$ is asymptotically
normal, with zero mean and unit variance. Consequently
$$
\ln W_{\mu} =\ln\mu +O_p(1).\tag 4.3.17
$$
Let $y=y(N)\to\infty$ so slow that $s=\mu e^y=o(\ln N)$ as well. Using the right hand 
side of (4.3.16), we obtain
$$
\align
P_i\{\ln W_{N,\mu}\ge \ln\mu+y\}=&\,P_i\{W_{N,\mu}\ge\mu e^y\}\\
=&\,P\{W_{\mu}\ge \mu e^y\}+O\bigl((s+\mu)/\ln N\bigr)\\
=&\,P\{\ln W_{\mu}\ge \ln\mu+y\}+O\bigl((s+\mu)/\ln N\bigr)=o(1).\tag 4.3.18
\endalign
$$
Analogously, the left hand side of (4.3.16) delivers
$$
\lim_{N\to\infty}P_i\{\ln W_{N,\mu}\ge \ln\mu-y\}=0.\tag 4.3.19
$$
The relations (4.3.17)-(4.3.19), together with (4.3.7) prove (4.3.10). 
\qed
\bi
The proof of Lemma 4.3.3 and Corollary 4.3.4 is similar and we omit it.
\bi

Recall that $N\in [N_1,N_2]$, $N_{1,2}=\lfloor \alpha\nu\pm 2\nu^{1/2}\ln \nu
\rfloor$.
Introduce $N_0=\lfloor\alpha\nu\rfloor$. It is easy to check that the
proof of Lemma 4.3.1 and Corollary 4.3.2 goes through with very minor changes
if, instead of (4.3.5), we define an integer $n$ by
$$
n=\frac{\ln\bigl[\lambda^{-1}BN_0(\ln N_0/\ln z_*)^2\bigr]}{\ln (1/z_*)}.
$$
(The key is that
$$
N(\ln N)^2=\bigl(1+O(N_0^{-1/2}\ln N_0)\bigr)N_0(\ln N_0)^2,
$$
uniformly for $N$ in question.) The same change can be made in the
formulation of Lemma 4.3.3 and Corollary 4.3.4 for the C-chain. This observation
coupled with the fact that $\hat X_+^{(\mu)}$, the $\mu$-th largest value of 
$\{X(t)\}_{0<t<\hat M}$, is sandwiched between those for $\{X(t)\}_{0<t\le N_1}$
and $\{X(t)\}_{0<t\le N_2}$, show that in Lemma 4.3.1, Corollary 4.3.2,
Lemma 4.3.3 and Corollary 4.3.4 we can put $\hat X_+^{(\mu)}$ instead of $X^{(\mu)}$. 
Below, by the relations (4.3.6) and (4.3.12) we will mean their modifications,
i.e. with $\hat X_+^{(\mu)}$ on their LHS.
\si

Turn to $\hat X^{(\mu)}$, the $\mu$-th largest value among
$X(0), X(1),\dots X(\hat M)$.  $X(0)$ has the distribution $\pi_1$ given
by either by (4.1.4) or by (4.1.14). Hence  
$$
P\{X(0)\ge n\}=\left\{\alignedat2
&O(nz_*^n),\quad&&\text{for cca},\\
&O(z_*^n),\quad&&\text{for C}.\endalignedat\right.\tag 4.3.20
$$
Now
$$
\hat X_+^{(\mu)}\le \hat X^{(\mu)}\le X(0) + \hat X_+^{(\mu)};
$$
so, for the cca case, we use $n$ defined by (4.3.5) and add an extra error term 
coming from (4.3.20), i.e.
$$
nz_*^n =O\left(z_*^{\frac{\ln(\nu\ln^2\nu)}{\lambda\ln(1/z_*)}}\ln\nu \right)=O(\nu^{-1}),
$$ 
to the RHS of (4.3.6), to obtain the corresponding claim for $\hat X^{(\mu)}$.
Likewise, in the C-case we need to add an error term $O(\nu^{-1/2})$ to the RHS
of (4.3.12). Again, we will refer to these new relations as (4.3.6) and (4.3.12).

But then, according to Theorem 4.2.1, the
$\mu$-th largest among the {\it parts\/} $Y_1,\mathbreak\dots,Y_{\hat\Cal M}$ of
the corresponding random composition can replace $\hat X^{(\mu)}$ on the
LHS of (4.3.6) and (4.3.12) respectively. These are our newest (4.3.6) and
(4.3.12).

Finally, if we include the 
rightmost parts $Y_{\hat\Cal M+1},Y_{\hat\Cal M+2},\dots$, it will not 
substantially affect the
the limiting behavior of the $\mu$-th largest overall part
either. Here is why. The number of these parts is $m:=\lceil\ln^2\nu\rceil$, at most. 
The total number of parts is q.s. of order $\nu\gg m$, which means the last $m$ parts
are q.s. well defined. Those parts, read from right to
left, and the first $m$ parts, read from left to right, are equidistributed.
By Theorem 4.2.1, these $m$ first parts are within the total variation
distance $O(\nu^{-K})$, ($\forall\,K>0$), from $Z(1),\dots, Z(m)$. We know
that $Z(1)$ has the distribution $\pi_1$. Since
$$
\sup_{i\in\Bbb N,\,t\ge 1}\sum_{k\ge n}p^t(i,k)\le 
\sup_{i\in\Bbb N}\sum_{k\ge n}p(i,k),\quad t\ge 1,
$$
we see that
$$
P\{Z(t)\ge n\}\le \sup_{i\in\Bbb N}\sum_{k\ge n}p(i,k),\quad t\ge 2.
$$
In view of (4.1.4) and (4.2.3), we obtain then: for the cca chain,
$$
\align
P\left\{\max_{1\le t\le m}Z(t)\ge\frac{\ln\bigl[\lambda^{-1}BN_0(\ln N_0/\ln z_*)^2
\bigr]}{\ln(1/z_*)}\right\}
=&\,O\left(\bigl(\ln^4\nu\bigr) z_*^{\frac{\ln(\nu\ln^2\nu)}{\lambda\ln(1/z_*)}}\right)\\
=&O\bigl(\nu^{-1}\ln^3\nu\bigr).
\endalign
$$
For the C-chain, the analogous probability is of order $\nu^{-1/2}\ln^2\nu$.
Therefore, by adding yet another error terms $O(\nu^{-1}\ln^3\nu)$ and $O(\nu^{-1/2}
\ln^2\nu)$ to the RHS of (4.3.6) and (4.3.12) (where $N=N_0=\lfloor\alpha\nu\rfloor$,
of course), we obtain the limiting distributions of the $\mu$-th largest part 
of both random compositions, together with explicit error terms. (For the
cca composition, the order of the total error term remains unchanged, i.e. 
$O\bigl((\lambda+\mu)/\ln\nu\bigr)$.
\bi

In summary, we have proved the following. 
\proclaim{Theorem 4.3.5} For a random composition $\bold Y$
of $\nu$, let $Y^{(\mu)}$ denote the $\mu$-th largest part. Let $N_0=\lfloor
\alpha\nu\rfloor$, $\alpha$ being defined in (4.1.17). Let $W_{\mu}$
be the sum of $\mu$ independent exponentials with unit mean. 
{\bf (i)\/} For the random cca composition, define $W_{\nu,\mu}$
by
$$
Y^{(\mu)}=\frac{\ln\bigl[W_{\nu,\mu}^{-1}BN_0(\ln N_0/\ln z_*)^2\bigr]}
{\ln(1/z_*)},
$$
$B$ being defined in (4.3.3). Then, for $s=o(\ln \nu)$ such that
$$
\frac{\ln\bigl[s^{-1}BN_0(\ln N_0/\ln z_*)^2\bigr]}
{\ln(1/z_*)}\in \Bbb N,
$$
we have 
$$
P\{W_{\nu,\mu}\ge s\}=P\{W_{\mu}\ge s\}+O\bigl[(s+\mu)/\ln\nu\bigr].
$$
\noindent
{\bf (ii)\/} For the random C-composition, define $W_{\nu,\mu}$ by
$$
Y^{(\mu)}=\frac{\ln\bigl[W_{\nu,\mu}^{-1}BN_0\bigr]}{\ln(1/z_*)},
$$
$B$ being defined in (4.3.4). Then, for $s=o(\nu^{1/2})$ such that
$$
\frac{\ln\bigl(s^{-1}BN_0\bigr)}{\ln(1/z_*)}\in \Bbb N,
$$
we have
$$
P\{W_{\nu,\mu}\ge s\}=P\{W_{\mu}\ge s\}+O\bigl(\nu^{-1/2}\ln^2\nu+(s^2+\mu^2)/\nu).
$$
\endproclaim
Here is a cruder estimate implied by Theorem 4.3.5.
\proclaim{Corollary 4.3.6} {\bf (i)\/} For the random cca composition,
$$
Y^{(\mu)}=\frac{\ln\bigl(\mu^{-1}\nu\ln^2\nu\bigr)}{\ln(1/z_*)}+O_p(1),\quad
(\mu=o(\ln\nu)).
$$
\noindent
{\bf (ii)\/} For the random C-composition,
$$
Y^{(\mu)}=\frac{\ln\bigl(\mu^{-1}\nu\bigr)}{\ln(1/z_*)}+O_p(1),\quad 
(\mu=o(\nu^{1/2})). 
$$
{\bf (iii)\/} So, for both cases,
$$
Y^{(1)}-Y^{(\mu)}=\frac{\ln\mu}{\ln(1/z_*)} +O_p(1),
$$
if $\mu=o(\ln\nu)$ and $\mu=o(\nu^{1/2})$ respectively.
\endproclaim
\bi
\bi
{\bf Acknowledgement.\/} It is my genuine pleasure to thank Guy Louchard
for introducing me to the random compositions and to his conjecture on
hitting times for Markov chains with uniformly exponential tails of the
row distributions. His strong belief in the conjecture, kind encouragement 
and insightful feedback helped 
to sustain my effort during months of work on this study. I thank a referee
for a painstaking effort to evaluate the paper and for a series of
penetrating critical comments.
\bi

\bi
\si
{\bf Appendix.\/}
\si
{\bf Proof of Lemma 4.1.1.\/} Consider the case of the random C-composition.
The BGF of $T(\mu,\nu)$, the number of C-compositions of $\nu$ with $\mu$ parts, and
$\nu$ and $\mu$ marked by $z$ and $w$ respectively, is given by (4.1.9)-(4.1.10):
$$
f(w,z)=-1+\frac{1}{h(w,z)},\quad h(w,z)=1-\sum_{j\ge 1}\frac{wz^j}{1+wz^j}.
$$
This bivariate series converges for $|z|<1$ and $|w|<1/|z|$. So, choosing $r_1<1$, and 
$r_2<1/r_1$, we have
$$\align
P\{\Cal M=\mu\}=&\,\frac{[z^{\nu}w^{\mu}]\,f(w,z)}{T(\nu)}\\
=&\frac{1}{T(\nu)}
\frac{1}{(2\pi i)^2}\oint\limits_{z\in \Cal C_1}\oint\limits_{w\in \Cal C_2}\frac{f(w,z)}
{z^{\nu+1}w^{\mu+1}}\,dwdz,
\endalign
$$
where $\Cal C_1$, $\Cal C_2$ are circles of radius $r_1$ and $r_2$ respectively.  In essence,
it is this formula that, via Bender's method [5], enabled Louchard [19] and
Louchard and Prodinger [21] to establish a sharp local limit theorem for
$\Cal M$ for the cca composition and the C-composition. Since our goal is to
bound the probability of large deviations, we use a considerably less analytical
argument, which is a bivariate extension of Chernoff's method. 

As a preparation, we need to define a differentiable extension of $z_*=0.57\dots$, the 
smallest-module root of $h(1,z)=0$. To this end, we compute
$$
h_z(1,z)=-\sum_{j\ge 1}\frac{jz^{j-1}}{(1+z^j)^2},\quad h_w(1,z)=-
\sum_{j\ge 1}\frac{z^j}{(1+z^j)^2}.\tag A.1
$$
So $h_z(1,z)<0$, $h_w(1,z)<0$ for $z\in (0,1)$. By continuity of $h_z(z,w)$, 
$h_w(z,w)$, we obtain: there exists $\varepsilon \in (0,1-z_*)$ such that (1) 
$(z_*+\varepsilon)(1+\varepsilon)<1$, and (2)
$$
h_z(z,w)<0,\,\,h_w(z,w) <0,\quad \forall\,(z,w)\in \Bbb R^2_+:\, z\le z_*+\varepsilon,\,
|w-1|\le\varepsilon.\tag A.2
$$
Consequently, for $|w-1|\le \varepsilon$, the equation $h(z,w)=0$ has a unique 
root $z=z(w)$, of multiplicity $1$, in
$[0,z_*+\varepsilon]$, which is infinitely differentiable as a function of $w$, and 
$z(1)=z_*$.
In particular,
$$
z^\prime(w)=-\frac{h_w(z(w),w)}{h_z(z(w),w)}<0,
$$
that is $z(w)$ is strictly decreasing. So $z(w)>z_*$ for $w<1$, and $z(w)<z_*$ for 
$w>1$.
\bi

Now, the series for the bivariate
generating function $f(w,z)$ converges for $|w-1|\le\varepsilon$ and $|z|<z(w)$. Since all the 
coefficients in the series  are nonnegative,
$$
\sum_{\ell\ge m}[z^\nu w^{\ell}]\,f(w,z)\le \frac{f(w,z)}{z^{\nu} w^m},\quad w\in [1,1+\varepsilon_0],\,
z\in (0, z(w)).
$$
Likewise
$$
\sum_{\ell\le m}[z^{\nu}w^{\ell}]f(w,z)\le \frac{f(w,z)}{z^{\nu} w^m},\quad w\in [1-\varepsilon_0,1],\,
z\in (0, z(w)).
$$
Here, by the definition of $f(w,z)$ and $z(w)$,
$$
f(w,z)\le\frac{c}{z(w)-z},\quad z<z(w).
$$
Therefore, for each $m$,
$$
P(\Cal M\ge m)\le c\frac{z^{-n}w^{-m}(z(w)-z)^{-1}}{T(\nu)},\quad w\in [1,1+\varepsilon_0],\,\,
z\in (0, z(w)), \tag A.3
$$
and
$$
P(\Cal M\le m)\le c\frac{z^{-n}w^{-m}(z(w)-z)^{-1}}{T(\nu)},\quad w\in [1-\varepsilon_0,1],\,\,
z\in (0, z(w)).\tag A.4
$$
Consider (A.3). To get the most out of this upper bound we need to determine $z$ and $w$ that minimize
the RHS, i. e.
$$
H^{(m)}(w,z):=-\nu\ln z-m\ln w-\ln (z(w)-z).
$$
Let us find a stationary point $(\bar w,\bar z)$ of $H^{(m)}(w,z)$ in the region 
$w\in [1,1+\varepsilon_0]$, $z\in (0,z(w))$. From the equations
$$
\align
H^{(m)}_z&=-\frac{\nu}{z}+\frac{1}{z(w)-z}=0,\\
H^{(m)}_w&=-\frac{m}{w}-\frac{z^\prime(w)}{z(w)-z}=0,
\endalign
$$
we obtain that 
$$
\bar z=\frac{\nu}{\nu+1}z(\bar w),
$$
where $\bar w=\bar w(m)$ must be a root of
$$
\frac{wz^\prime(w)}{z(w)}=-\frac{m}{\nu+1}.\tag A.5
$$
The equation (A.5) has a solution $w=1$ if 
$$
m=\bar m:=(\nu+1)\,\mu,\quad \mu:=-\frac{z^\prime(1)}{z(1)}.
$$
Furthermore, in [19] it was shown that
$$
\left.\frac{d}{dw}\frac{wz^\prime(w)}{z(w)}\right|_{w=1}=\frac{z^{\prime\prime}(1)}{z(1)}
-\mu-\mu^2
$$
is negative; this is $-\beta$, $\beta$ defined in (4.1.16). 
Since
$$
\frac{d}{dm}\left(-\frac{m}{\nu+1}\right)=-\frac{1}{\nu+1}<0
$$
as well, for 
$$
0\le m-\bar m=o(\nu),
$$
the equation (A.5) defines a strictly increasing $\bar w(m)$; so $\bar w(m)>1$ for
$m>\bar m$. More precisely
$$
\align
\bar w(m)&=1+\frac{\beta}{\nu+1}(m-\bar m)+O((m-\bar m)^2/\nu^2)\\
&=1+\frac{\beta}{\nu}(m-\mu\,\nu)+O((m-\mu\,\nu)^2/\nu^2).
\endalign
$$
Now
$$
\align
H^{(\bar m)}(\bar w(\bar m),\bar z(\bar w(\bar m)))&=-\nu\ln\left(\frac{\nu}
{\nu+1}z_*\right)-
\ln\left(\frac{z_*}{\nu+1}\right)\\
=&-\nu\ln z_*+\ln \nu+O(1).
\endalign
$$
Also
$$
\align
&\frac{d}{dm}H^{(m)}(\bar w(m),\bar z(m))\\
=&H^{(m)}_m(\bar w(m),\bar z(m))+
H^{(m)}_w(\bar w(m),\bar z(m))+H^{(m)}_z(\bar w(m),\bar z(m))\\
=&H^{(m)}_m(\bar w(m),\bar z(m))=-\ln \bar w(m),
\endalign
$$
which implies that
$$ 
\left.\frac{d}{dm}H^{(m)}(\bar w(m),\bar z(m))\right|_{m=\bar m}=-\ln\bar w(\bar m)=0,
$$
and also that
$$
\align
\left.\frac{d^2}{dm^2}H^{(m)}(\bar w(m),\bar z(m))\right|_{m=\bar m}=&
\left.-\frac{\bar w^\prime(m)}{\bar w(m)}\right|_{m=\bar m}\\
=&-\bar w^\prime(\bar m)\\
=&-\frac{\beta}{\nu+1}.
\endalign
$$
Therefore, for $0\le m-\bar m=o(\nu)$, 
$$
\align
H^{(m)}(\bar w(m),\bar z(m))=&-\nu\ln z_*+\ln \nu-(1+o(1))\frac{\beta}{2(\nu+1)}
(m-\bar m)^2+O(1)\\
\le &-\nu\ln z_* +\ln \nu-\frac{\beta}{3\nu}
(m-\mu\,\nu)^2+O(1).
\endalign
$$
Using this bound in (A.3) for $w=\bar w(m)$, $z=\bar z(\bar w(m))$, and recalling that
$T(\nu)$ is of order $z_*^{-\nu}$, we obtain:
$$
P(\Cal M\ge m)\le c \nu\exp\left(-\frac{\beta}{3\nu}(m-\mu\,\nu)^2\right),
\quad 0<m-\mu\, \nu=o(\nu).
$$
Likewise
$$
P(\Cal M\le m)\le c \nu\exp\left(-\frac{\beta}{3\nu}(m-\mu\,\nu)^2\right),
\quad 0<\mu\,\nu-m=o(\nu).
$$

The case of the random cca composition is quite analogous, so we omit the proof.
\qed
\enddocument